\documentclass[12pt,a4paper,leqno]{amsart}
\usepackage{amsmath}
\usepackage{amssymb}

\newtheorem{theorem}{Theorem}[section]

\newtheorem*{theoremD}{Theorem  D} 
\newtheorem{lemma}[theorem]{Lemma}
\newtheorem{proposition}[theorem]{Proposition}
\newtheorem{corollary}[theorem]{Corollary}
\newtheorem{remark}[theorem]{Remark}
\newtheorem{definition}[theorem]{Definition}

\newcommand{\C}{{\bf C}}       
\newcommand{\B}{{\bf B}}       
\newcommand{\R}{{\bf R}}       
\newcommand{\N}{{\bf Z}_{+}}   
\newcommand{\Z}{{\bf Z}}       

\catcode`\@=11 \@addtoreset{equation}{section} \catcode`\@=12

\title[Carleson measures for weighted holomorphic Besov spaces in $\C^n$] { Carleson measures  for weighted   holomorphic Besov spaces }
\author[Cascante]{Carme Cascante}
\address{Departament de Matem\`atica Aplicada i An\`alisi,
Facultat de Matem\`atiques, Universitat de Barcelona, Gran Via
585, 08071~Barcelona, Spain} \email{cascante@ub.edu}

\author[Ortega]{Joaquin M. Ortega }
\address{ Departament de Matem\`atica Aplicada i An\`alisi,
Facultat de Matem\`atiques,  Universitat de Barcelona, Gran Via
585, 08071~Barcelona, Spain} \email{ortega@ub.edu}

\thanks{Both authors partially supported by DGICYT Grant
MTM2005-08984-C02-02, and DURSI Grant 2005SGR 00611
.}

\date{}

\begin{document}
\begin{abstract}
We obtain characterizations of positive Borel measures $\mu$ on
$\B^n$
 so that some weighted  holomorphic Besov spaces $B_s^p(w)$ are imbedded in  
  $L^p(d\mu)$,  where $w$ is a  $B_p$ weight in  the unit ball of $\C^n$. 
\end{abstract}
\maketitle

\section{Introduction}  \label{introduction}


If $w$ is a weight in $\B^n$, the unit ball in $\C^n$, and $p>0$, $s\in\R$,  the space $B_s^p(w)$ consists of holomorphic functions on $\B^n$ for which $$||f||_{B_s^p(w)}^p=\int_{\B^n}|(I+R)^kf(y)|^p(1-|y|^2)^{(k-s)p-1} w(y)dv(y)<+\infty, $$
 for some $k\in \Z_+$, $k>s$. Here $dv$ is the normalized Lebesgue measure on $\B^n$.  As it happens in the unweighted case, it can be shown  that for adequate weights if the above integral is finite for some $k>s$, then it is also finite for any $k>s$ (see section $3$).
 
In this paper we consider Carleson measures for  weighted holomorphic Besov space  $B_s^p(w)$, that is, the positive Borel measures $\mu$ on $\B^n$, the   unit ball in $\C^n$, for which the weighted holomorphic Besov space  space $B_s^p(w)$ is imbedded in $L^p(d\mu)$. 
 
 For some particular cases of weighted Besov spaces the characterization of the corresponding Carleson measures is known. For instance, when $s<0$, no derivative is necessarily involved in the definition of the norm of $B_s^p(w)$, and it is in fact a weighted Bergman space. In that case, if   $w(z)=(1-|z|)^\alpha$,  where $\alpha -sp>0$,  $\mu$ is a Carleson measure for $B_{s}^p(w)$ if and only if there exists $C>0$ such that for any $\eta\in{\bf S}^n$, and $R>0$, $\mu(T(\eta,R))\leq C R^{n+ \alpha -sp}$, where $T(\eta, R)=\{ z\in \B^n,\,\, |1-z\overline{\eta}| <R\}$  (see \cite{oleinik-pavlov}, \cite{stegenga} in dimension $1$, and \cite{luecking} in dimension $n>1$ among others). This result can be extended to $\alpha=sp$ and $p\leq 2$ (see for instance the survey \cite{zhaozhu}). On the other hand, if $n+\alpha-sp<0$, it is well known that the space $B_s^p(w)$ consists of regular functions, and the Carleson measures in those cases are just the finite ones.

 Let's finally mention that if $n=1$, \cite{arcozzirochbergsawyer}  have studied the Carleson measures for  $B_{\frac{1}{p'}}^p(w)$, where $1<p$, $0<s<1$, and $w$ is a weight in $\B^1$ in the class   $B_p$   (see \cite{bekolle} for a definition)  satisfying some additional regularity conditions on $w$. 
  
  The main purpose of these paper is to obtain characterizations of Carleson measures for weighted Besov spaces in dimension $n>1$.  We remark that they also appear in a natural way in the study of Carleson measures for Hardy-Sobolev spaces on general domains and such that the support of the measure is included on a manifold. For instance, if $D=\{ (z,y)\,;\, z,y\in \B^1,\, |y|\leq \varphi(z)\}$, where $\varphi$ is a nonnegative function in ${\mathcal C}^1(\B^1)$, the study of Carleson measures for $H_s^p(D)$ supported on $y=0$, $s>\frac1{p}$ leads to the study of Carleson measures on $B_{s-\frac1{p}}(w)$ where $w={\displaystyle \varphi(4\frac{\partial \varphi}{\partial z}\frac{\partial \varphi}{\partial \overline{z}}+1)^\frac12}$. Although we will not  study them in detail, our techniques are inspired in these facts. 
  
We observe that when $w\equiv 1$,  the holomorphic Besov space $B_s^p$ can be viewed as a restriction to $\B^n$ of the Hardy-Sobolev space $H_{s+\frac1{p}}^p(\B^{n+1})$ (see \cite{ortegafabrega}). This fact allows to reduce the study  of Carleson measures for $B_s^p$ to   the study of Carleson measures for Hardy Sobolev spaces in one more dimension. When $n-sp<1$, there exists characterizations of Carleson measures for $H_{s+\frac1{p}}^p(\B^{n+1})$ ( see \cite{cascanteortega2})  which give characterizations for Carleson measures for $B_s^p(\B^n)$. Other authors (see \cite{arcozzirochbergsawyer1} and references therein) have also obtained other type of characterizations for $H_s^2(\B^n)$ when   $n-2s\leq 1$.

 Let's recall the main facts about Carleson measures on weighted Hardy-Sobolev spaces that we will need. The weighted Hardy-Sobolev space $H_s^p(w,\B^n)$,   $0\leq s, p<+\infty$, consists of those functions $f$ holomorphic in $\B^n$ such that if ${\displaystyle f(z)=\sum_k f_k(z)}$ is its homogeneous polynomial expansion, and   ${\displaystyle (I+R)^s f(z)=\sum_k (1+k)^s f_k(z)},$  we have that
$${\displaystyle ||f||_{H_s^p(w)}=\sup_{0<r<1} ||(I+R)^s f(r\zeta)||_{L^p(w)}<+\infty}.$$
We also recall that a $w$ is in $A_p({\bf S}^n)$, $1<p<+\infty$,  if there exists $C>0$ such that for any nonisotropic ball $B\subset {\bf S}^n$, $B=B(\zeta,r)=\{ \eta\in{\bf S}^n\,;\, |1-\zeta\overline{\eta}|<r\,\}$,
$$\left( \frac{1}{|B|}\int_B w d\sigma\right)\left( \frac{1}{|B|}\int_B w^{\frac{-1}{p-1}}d\sigma \right)^{p-1}\leq C,$$
where $\sigma$ is the Lebesgue measure on ${\bf S}^n$ and $|B|$ denotes the Lebesgue measure of the ball $B$.

A weight $w$ in ${\bf S}^n$ is in  $D_\tau({\bf S}^n)$, if   there exists $C>0$ such that for any nonisotropic ball $B$ in ${\bf S}^n$, $w(2^kB)\leq C 2^{k\tau}w(B)$. Analogously to what it happens with weights in $\R^n$, the fact that a weight is in $A_p({\bf S}^n)$ implies that it is in   $D_\tau({\bf S}^n)$
 for    $\tau= np$. 

We denote by $K_s$ the nonisotropic potential operator defined by
$$K_s[f](z)=\int_{{\bf S}^n} \frac{f(\eta)}{|1-z\overline{\eta}|^{n-s}}d\sigma(\eta), \,\, z\in\overline{\B}^n.$$
It has been shown in \cite{cascanteortega3}
 \begin{theorem}[\cite{cascanteortega3}]\label{weightedtraceinequality}
 Let $1< p<+\infty$,  $w$  an $A_p$-weight, and $\mu$ a finite positive Borel measure on $\B^n$. Assume that $w$ is in $D_\tau$ for some $0\leq \tau-sp<1$. We then have that the following statements are equivalent:
 
 (i) $||f||_{L^p(d\mu)}\leq C|| f||_{H_s^p(w\B^n)}$.
 
 (ii) $||K_\alpha(f)||_{L^p(d\mu)}\leq C||f||_{ L^p(w)}$.

\end{theorem}
 
  We remark that this equivalence is quite useful in many applications, since allows to work with a positive kernel.

Let's now state the main result in this paper. Recall that a weight $w$ is in $B_p(\B^n)$,   if there exists $C>0$ such that for any tent $T(\zeta,R)$, $\zeta\in{\bf S}^n$,
 $$\frac1{v(T(\zeta,R))}\int_{T(\zeta,R)}\omega dv \left( \frac1{v(T(\zeta,R))}\int_{T(\zeta,R)}\omega^{-(p'-1)}dv\right)^{p-1} \leq C.$$ 
 We
 introduce a pseudodistance in $\overline{\B}^n$
defined by 
$$\rho(z,w)=|1-z\overline{w}|-\sqrt{1-|z|^2} \sqrt{1-|w|^2}.$$

 \begin{theoremD}
 Let $1<p<+\infty$, $w$ a $B_p(\B^n)$ weight. Assume that $w$ satisfies a doubling condition of order $\tau+1$, $\tau<1+sp$, for the pseudodistance $\rho$. Let $\mu$ be a  positive Borel measure on $\B^n$. We then have that the following statements are equivalent:

(i) There exists $C>0$ such that for any $f\in B_{s}^{p}(w)$,
$$||f||_{L^p(d\mu)}\leq C||f||_{B_{s}^{p}(w)}.$$

(ii) There exists $C>0$ such that for any $f\in L^p(wdv)$, 
$$||\int_{\B^n} \frac{f(y)dv(y)}{(1-z\overline{y})^{n+1-(s+\frac1{p})}}||_{L^p(d\mu)}\leq C||f||_{L^p(wdv)}.$$

(iii) There exists $C>0$ such that for any $f\in L^p(wdv)$, 
$$||\int_{\B^n} \frac{f(y)dv(y)}{|1-z\overline{y}|^{n+1-(s+\frac1{p})}}||_{L^p(d\mu)}\leq C||f||_{L^p(wdv)}.$$
 \end{theoremD}
 
 The paper is organized as follows: In section 2 we introduce the class of weights we will consider. We obtain all the properties on the weights  needed in the proof of Theorem D. In section 3 we study the general properties of the weighted Besov spaces $B_s^p(w)$ and in section 4 we will give the proof of theorem D.

 Finally, the usual remark on  notation: we will adopt the  convention
of using the same letter for various absolute constants whose
values may change in each occurrence, and we will write $A\preceq
B$ if there exists an absolute constant $M$ such that $A\leq MB$.
 We will say that two quantities $A$ and $B$ are equivalent if both
$A\preceq B$ and $B\preceq A$, and, in that case, we will write
$A\simeq B$.

    

\section{weights in $\B^n$}
\label{section2}

 Our approach to the study of weighted Besov spaces in $\B^n$ uses their immersion in holomorphic spaces defined in ${\B^{n+1}}$ via the natural projection $\Pi:{\B}^{n+1}\rightarrow  \B^n$, given by $\Pi(z_1,\cdots,z_{n+1})=(z_1,\cdots,z_n)$. It is then convenient to consider a pseudodistance  in ${\overline{\B}}^n$  deduced from the hyperbolic pseudodistance in ${\bf S}^{n+1}$.

If $z,w\in\overline{\B}^n$, let 
\begin{equation*}\begin{split}&\rho(z,w)=\inf_{\varphi,\,\theta\in [0,2\pi)}|1-z\overline{w}-\sqrt{1-|z|^2}e^{i\varphi}\sqrt{1-|w|^2}e^{i\theta}|\\&=\inf_{\theta\in [0,2\pi)}|1-z\overline{w}-\sqrt{1-|z|^2}\sqrt{1-|w|^2}e^{i\theta}|=|1-z\overline{w}|-\sqrt{1-|z|^2}\sqrt{1-|w|^2}.\end{split}\end{equation*} 

Observe that $\rho(z,w)$ is just  the infimum of the Korany pseudodistances of the antiimages by the mapping $\Pi$ of the points $z,w$.

 Let us see a suggestive expression of $\rho$ that has $|1-z\overline{w}|$ as a factor. Let   $P_a$ be the orthogonal projection of $\C^n$ onto the subspace $[a]$ generated by $a$ and $Q_a=Id-P_a$   the projection onto the orthogonal complement of $[a]$. If $a\in\B^n$, and $\varphi_a$ is the automorphism in $\B^n$ which interchanges $a$ and $0$, given by
$$
\varphi_a(z)=\frac{a-P_a(z)-(1-|a|^2)^\frac12 Q_a(z)}{1-z\overline{a}},
$$
then (see for instance theorem 2.2.2 in \cite{rudin})
$$1-|\varphi_a(z)|^2= \frac{(1-|a|^2)(1-|z|^2)}{|1-z\overline{a}|^2}.$$
This fact gives  that $\rho(z,w)\simeq |1-z\overline{w}||\varphi_z(w)|^2$. Indeed,
\begin{align*}&\rho(z,w)=|1-z\overline{w}|-\sqrt{1-|z|^2}\sqrt{1-|w|^2}\\&=|1-z\overline{w}|\left(1-\frac{\sqrt{1-|z|^2}\sqrt{1-|w|^2}}{|1-z\overline{w}|}\right)\simeq |1-z\overline{w}||\varphi_z(w)|^2.\end{align*}
 
 In the  following lemma we show that $\rho$ is a pseudodistance whose balls $U(z,R)$  are ''equivalent'', in a sense that we will precise, to polydisks  of size $ R+R^\frac12 (1-|z|^2)^\frac12$ in the complex normal direction and of size $R^\frac12$ in the complex-tangential directions. In order to distinguish the Lebesgue measure in ${\bf S}^{n+1}$ and in $\B^n$, we will write $v(E)$ the volume measure of a measurable subset $E$ in $\overline{\B}^n$, whereas $|F|$ will stand for the Lebesgue measure of a measurable subset $F\subset{\bf S}^{n+1}$.
 
\begin{lemma}\label{pseudodistance}
\begin{itemize}\item[(i)] $\rho$ is a pseudodistance in $\overline{\B}^n$.

\item [(ii)] Let $z\in \overline{\B}^n$, $0<R<1$, and $U(z,R)=\{ w\in \overline{\B}^n\,;\, \rho(z,w)<R\}$. Let $P(z,R)$ be the polydisk in $\overline{\B}^n$ centered at $z$, of size $ R+R^\frac12 (1-|z|^2)^\frac12$ in the complex normal direction and of size $R^\frac12$ in the complex-tangential directions.  Then there exists $C>0$ such that  $P(z,\frac{R}{C})\subset U(z,R)\subset P(z,CR)$. In particular, $v(U(z,R))\simeq R^{n}(R+(1-|z|^2))$.
\end{itemize}
\end{lemma}
{\bf Proof of lemma \ref{pseudodistance}:}\par
 
 Let $\theta \in [0,2\pi)$. If $z\in\B^n$, let $\Pi_\theta^{-1}(z)=(z,\sqrt{1-|z|^2}e^{i\theta})$. We then have that for any $\theta_0\in[0,2\pi)$,  $\rho(z,w)=\inf_{\theta}|1-\Pi^{-1}_\theta(z)\overline{\Pi^{-1}_{\theta_0}(w)}|$, expression from which we easily obtain that $\rho$ is a pseudodistance. That gives (i).
 
 Let us prove (ii). Let $z\in\overline{\B}^n$, and $R>0$. A unitary change of variables gives that, without loss of generality,  we may assume that $z=(r,\stackrel{n-1)}{0,\cdots,0})$, $0\leq r\leq 1$.  We begin showing that $P(z,R)\subset U(z,CR)$, for some fixed constant $C>0$. Let us consider first the case $R\leq 1-r^2$. If $w=(w_1,\cdots,w_n)\in P(z,R)$, the definition of the polydisk gives that $|r-w_1|\leq R+R^\frac12 (1-r^2)^\frac12$, and $|w_i|\leq R^\frac12$, $i=2,\cdots,n$. Then
 \begin{equation*}\begin{split}&\rho(z,w)=|1-r\overline{w_1}|-\sqrt{1-r^2}\sqrt{1-|w|^2}=\frac{|r-{w_1}|^2+(|w_2|^2+\cdots+|w_n|^2)(1-r^2)}{ |1-r\overline{w_1}|+\sqrt{1-|w|^2}\sqrt{1-r^2}}  \\&\simeq \frac{|r-{w_1}|^2+(|w_2|^2+\cdots+|w_n|^2)(1-r^2)}{ |1-r\overline{w_1}|}\preceq 
 \frac{R^2+R(1-r^2)+R(1-r^2)}{1-r^2}\preceq R.
 \end{split}\end{equation*}
 Assume now that $(1-r^2)\leq R$. We have that
 $$\rho(z,w)\leq |1-r\overline{w_1}|\simeq (1-r^2)+|r-w_1|\preceq 
 (1-r^2)+R+R^\frac12(1-r^2)^\frac12\preceq R.
 $$
Hence in any case we have shown that $P(z,R)\subset U(z,CR)$.

Conversely, let  $w\in\B^n$ such that $\rho(z,w)<R$. The previous argument gives that  
$$ \rho(z,w)\simeq \frac{|r- {w_1}|^2+(|w_2|^2+\cdots+|w_n|^2)(1-r^2)}{|1-r\overline{w_1}|}\preceq R.$$ In particular, \begin{equation}\label{desigu}\frac{|r- {w_1}|^2}{(1-r^2)+|r-w_1|}\preceq R.\end{equation}
If $|r-w_1|\leq (1-r^2)$, we have that $(1-r^2)+|r-w_1|\simeq (1-r^2)$, and we deduce from (\ref{desigu}) that    $|r-w_1|\preceq R^\frac12(1-r^2)^\frac12$.  If on the other hand, $(1-r^2)\leq |r-w_1|$, (\ref{desigu}) gives that $|r-w_1|\preceq R$. Thus in any   case we deduce that $|r-w_1|\preceq R+R^\frac12(1-r^2)^\frac12$. In order to finish we have to check that $|w_i|^2\preceq R$, $i=2,\cdots, n$. It is clear that this is the case if $|r-w_1|\leq (1-r^2)$, since then $\frac{(|w_2|^2+\cdots |w_n|^2)(1-r^2)}{1-r^2}\preceq R$. So we may assume that $(1-r^2)\leq |r-w_1|$. Since we have shown that in that case $|r-w_1|\preceq R$, we obtain
$$|w_2|^2+\cdots+|w_n|^2\leq 1-|w_1|^2\preceq |1-w_1|\preceq (1-r^2)+|r-w_1|\preceq |r-w_1|\preceq R.$$ 
The affirmation on the volume of the balls is obvious from the above.\qed

Observe that if $0<\varepsilon<1$, and $U_\varepsilon(z)=U(z,\varepsilon(1-|z|))$, we have that $U_\varepsilon(z)$ are contained and contain   ellipsoids $E(z)=\{w\in\B^n\,;\, |\varphi_z(w)|<\varepsilon'\}$, $\varepsilon'<1$, where $\varphi_z$ is the automorphism in $\B^n$ that interchanges $z$ and $0$. This can be checked as follows: if $w\in U_\varepsilon(z)$, $\rho(z,w)<\varepsilon(1-|z|)$, we  have that $|1-z\overline{w}||\varphi_z(w)|^2\preceq \varepsilon (1-|z|)$, and consequently $|\varphi_z(w)|\preceq \sqrt{\varepsilon}$. Conversely, if $|\varphi_z(w)|<\varepsilon^\frac12$, we have that 
$$ |\varphi_a(z)|^2= 1-\frac{(1-|a|^2)(1-|z|^2)}{|1-z\overline{a}|^2}<\varepsilon.$$
Thus,
$$\frac{|1-z\overline{a}|^2}{(1-|a|^2)(1-|z|^2)}<\frac1{1-\varepsilon},$$
and 
$$ |1-z\overline{a}|^2 |\varphi_a(z)|^2 <\frac{\varepsilon }{1-\varepsilon} (1-|a|^2)(1-|z|^2) .$$
Hence
$$ |1-z\overline{a}|  |\varphi_a(z)|^2 <\frac{\varepsilon }{1-\varepsilon}(1-|z|^2).$$

In particular, we deduce that provided $\varepsilon$ is small enough, there exists a constant $C>0$, such that if $w\in U_\varepsilon (z)$, we have that $\frac1{C}(1-|z|^2)\leq   (1-|w|^2)\leq C(1-|z|^2)$.

\begin{definition} \label{definitionAp}
 We say that a weight $w$ is in $A_p(\B^n)$, $1<p<+\infty$, if there exists $C>0$
 such that for any  ball $U=U(z,R)$ in $\overline{\B}^n$ associated to the pseudodistance $\rho$, 
$$\left( \frac{1}{v(U)}\int_U w dv\right)\left( \frac{1}{v(U)}\int_U w^{-(p'-1)}dv \right)^{p-1}\leq C.$$
  
 \end{definition}

 In the following lemma, we will obtain a characterization of  weights in $A_p(\B^{n})$ in terms of their ''lifting'' to ${\bf S}^{n+1}$. We recall that a weight $\eta$ in ${\bf S}^{n+1}$ is in $A_p({\bf S}^{n+1})$ if there exists $C>0$ such that for any nonisotropic ball $B(\zeta,R)=\{z\in{\bf S}^{n+1}\,;\, |1-\zeta\overline{z}|<R\}$, $\zeta\in{\bf S}^{n+1}$, $R>0$,
 $$\left( \frac{1}{|B(\zeta,R)|}\int_{B(\zeta,R)} \eta d\sigma\right)\left( \frac{1}{|B(\zeta,R)|}\int_{B(\zeta,R)} \eta^{-(p'-1)}d\sigma \right)^{p-1}\leq C.$$

 \begin{lemma}\label{lifted}
 Let $1<p<+\infty$, $n\geq 1$, and $w$ be a weight in $\B^n$. We then have that $w$ is an $A_p(\B^n)$  weight if and only if the lifted weight $w_l$ defined by $w_l(z_1,\cdots,z_{n+1})=w(z_1,\cdots,z_{n})$ is an $A_p$ weight in ${\bf S}^{n+1}$. 
 \end{lemma}
 {\bf Proof of lemma \ref{lifted}:}\par
 
 We begin proving that if $w$ is an $A_p(\B^n)$-weight, then $w_l$ is an $A^p$ weight in ${\bf S}^{n+1}$. We recall that we have denoted by $U$ the balls in $\B^n$ with respect to the pseudodistance $\rho$, and we will denote by $B(z,R)=\{y\in{\bf S}^{n+1}\,; |1-z\overline{y}|<R\,\}$,  the ball  in ${\bf S}^{n+1}$ of center $z=(z_1,\dots,z_{n+1})\in{\bf S}^{n+1}$ and radious $R$.

   We consider first the particular case where $z_{n+1}=0$, i.e.  the center of the ball $B(z,R)$ lies  on ${\overline\B}^n$. By a suitable change of variables   we may assume that $z=(1,{0,\cdots,0)}$.
 Then  $\Pi^{-1}(U((1,\stackrel{n-1}{0,\dots,0}),R))=\{ y\in{\bf S}^{n+1}\,;\, |1-y_1|<R\}=B(z,R)$, and consequentely,
$$\int_{B(z,R)}w_ld\sigma= \int_{\Pi^{-1}(U((1,\stackrel{n-1}{0,\dots,0}), R))}w_ld\sigma=  \int_{U((1,\stackrel{n-1}{0,\dots,0}), R)}wdv,
 $$
 and the same argument holds for  $w_l^{-(p'-1)}$. These estimates, together with the fact that $w\in A_p(\B^n)$ and $v(P(1,\stackrel{n-1}{0,\dots,0}),R)\simeq R^{n+1}$, give that
 $$\frac1{R^{n+1}}\int_{B(z,R) }w_ld\sigma \left( \frac1{R^{n+1}}\int_{B(z,R)}w_l^{-(p'-1)}d\sigma\right)^\frac1{p'-1} \leq C.$$
 In fact, this argument can be applied to nonisotropic balls $B(z,R)$ in ${\bf S}^{n+1}$ satisfying that $d(z,{\bf S}^n)\leq R$, where ${\bf S}^n$ is the boundary of ${\overline\B}^n$. We just have to observe that in this case the ball $B(z,R)$ is included in a nonisotropic ball in ${\bf S}^{n+1}$ whose center lies in ${\bf S}^n$ and whose radius is comparable to $R$. 
 
 So  we may assume that  $R\leq d(z,{\bf S}^n)$. 
 Without loss of generality we also may assume that $z=(r,\stackrel{n-1}{0\dots,0},\sqrt{1-r^2}e^{i\theta}))$, for some $0<r<1$, $\theta\in[0,2\pi)$. The fact that $R\leq d(z,{\bf S}^n)$ gives that $R\preceq 1-r^2$, and consequently that $v(U((\stackrel{n-1}{r,0\dots,0)},R))\simeq R ^n(1-r^2) $.

 If we denote by $T$  the unitary map in ${\C}^{n+1}$ given by
 \[T(y_1,\dots,y_{n+1})= \left( \begin{array}{ccccc}r& 0&0&\cdots&-e^{-i\theta}\sqrt{1-r^2}\\0&1&0&\cdots&0\\ 0&0&1&\cdots&0\\ \vdots&\vdots&\vdots&\vdots&\vdots\\
 \sqrt{1-r^2}e^{i\theta}& 0&0&\cdots&r\end{array}\right)\left(\begin{array}{c}y_1\\y_2\\\vdots\\\vdots\\y_{n+1}  \end{array}\right),\]
  we have that $T(1,\stackrel{n}{0,\dots,0})=z$, and consequentely, $T(B((1,\stackrel{n}{0,\dots,0}),R)) =B(z,R)$.
 Since $B((1,\stackrel{n}{0,\dots,0}),R)=\{ (y,y_{n+1})=(y_1,\dots,y_n ,\sqrt{1-|y |^2}e^{i\varphi})\,;\, 0\leq\varphi<2\pi,\,\, |1-y_1|<R,\},$ we obtain that
 \begin{equation*}\begin{split}&\Pi(B(z,R))=\Pi(T(B((1,\stackrel{n}{0,\dots,0}),R))\\&=\{  (ry_1-e^{i(\varphi-\theta)}\sqrt{1-r^2}\sqrt{1-|y|^2},y_2,\dots,y_n)\,;\,\\&
 0\leq \varphi<2\pi\,, |1-y_1|<R,\, |y_2|^2+\cdots|y_n|^2<1-|y_1|^2\}.
 \end{split}\end{equation*}
 
 This description of the set $\Pi(B(z,R))$  gives that $|ry_1-e^{i(\varphi-\theta)}\sqrt{1-r^2}\sqrt{1-|y|^2}-r|\preceq R+\sqrt{1-r^2}R^\frac12$, and $|y_k|\preceq R^\frac12$. Hence,  Lemma \ref{pseudodistance} gives that \begin{equation}\label{liftedweight}\Pi(B(z,R))\subset U((r,\stackrel{n-1}{0,\dots,0)}, CR).\end{equation}

On the other hand, let  $y=(y_1,\dots,y_n)\in U(( r,\stackrel{n-1}{0,\dots,0)}, CR)$,  and   let $ y^\tau=(y_1,\dots,y_n,\sqrt{1-|y|^2}e^{i\tau}) $ be a  point  in $\Pi^{-1}(y).$
 The definition of the pseudodistance $\rho$ gives that there exists a point in $\Pi^{-1}((r,\stackrel{n-1}{0,\dots,0}))$ at distance from the point $y^\tau$ less than or equal to $CR$.
This holds if and only if \begin{equation*}\begin{split}&CR\geq \min_{\varphi} |1-y_1r-\sqrt{1-|y|^2}\sqrt{1-r^2}e^{i(\tau- \varphi)}|=  \\&|1-y_1r|-\sqrt{1-|y|^2}\sqrt{1-r^2}=\rho(( r,\stackrel{n-1}{0\dots,0)},y).\end{split} \end{equation*}
  
 Summarizing, the projection of the points in ${{\bf S}^{n+1}}$ at distance less that $R$ from the set $\Pi^{-1}((r,\stackrel{n-1}{0,\dots,0}))$ is included in $U(( r,\stackrel{n-1}{0\dots,0)},CR)$, 
 and on the other hand, the ball $U(( r,\stackrel{n-1}{0\dots,0)},R)$ is included in the projection of the set of points in ${{\bf S}^{n+1}}$ at distance from $\Pi^{-1}((r,\stackrel{n-1}{0,\dots,0}))$ less than $CR$. 
 
 Next, this set of points of ${{\bf S}^{n+1}}$ at distance less that $R$ from the set $\Pi^{-1}((r,\stackrel{n-1}{0,\dots,0}))$ is included in a union of balls of radious $R$ in a number which is of the order of $[\frac{1-|z|^2}{R}]$, and includes the same number of disjoint balls of radious comparable to $R$. The integral of the lifted weight $w_l$ on each of this balls is equivalent.

Altogether we obtain that
\begin{equation*}\begin{split}&\frac{1-r^2}{R}\int_{B(( r,\stackrel{n-1}{0,\dots,0},\sqrt{1-r^2}e^{i\theta}),R)} w_ld\sigma\simeq \int_{\Pi^{-1} (U(( r,\stackrel{n-1}{0,\dots,0,)}, CR ))}w_ld\sigma= \\& \int_{ U((r,\stackrel{n-1}{0,\dots,0,}),CR)}wdv,
\end{split}\end{equation*}
and, since in the case we are now considering $v(U((r,\stackrel{n-1}{0,\dots,0)})\simeq R ^n(1-r^2) $, we have that
\begin{equation*}\begin{split}&
\frac1{R^{n+1}} \int_{B(( r,\stackrel{n-1}{0,\dots,0},\sqrt{1-r^2}e^{i\theta}),R)} w_ld\sigma\simeq \frac1{(1-|z|^2)  R^n}\int_{U( ( r,\stackrel{n-1}{0,\dots,0,)}, CR)}wdv \\& \simeq\frac1{v(U( ( r,\stackrel{n-1}{0,\dots,0,)}, CR))}\int_{U( (r,\stackrel{n-1}{0,\cdots,0,)}, CR)}wdv,
\end{split}\end{equation*}
with a simmilar estimate for $w_l^{-(p'-1)}$. Hence, we have proved that if $w\in A_p(\B^n)$, then $w_l$ is an $A_p$ weight in ${\bf S}^{n+1}$.

Let now $w$ be a weight in $\B^n$ satisfying that $w_l$ is an $A_p$ weight in ${\bf S}^{n+1}$. The argument we have used before shows that if $z\in \overline{\B^n}$, $R>0$ and $U(z,R)$ is a ball in $\B^n$, such that $1-|z|^2\leq R$, we can reduce ourselves to the case where the point $z$ is in ${\bf S}^n$. Then $U(z,R)$ is just a tent centered at a point $z$ in ${\bf S}^n$, and  $U(z,R)=B((z,0),R)$. Consequently, $w$ satisfies the $A_p(\B^n)$ condition for these class of balls.

If $R<1-|z|^2$, $v(U(z,R))\simeq R^n(1-|z|^2)$, and again the argument used before gives then that for any $\theta$,
\begin{equation*}\begin{split}&
\frac1{R^n(1-|z|^2)}\int_{U(z,R)}wdv\simeq \frac1{R^n(1-|z|^2)}\int_{\Pi^{-1}(U(z,R))}w_ld\sigma
\\&\simeq
\frac1{R^n(1-|z|^2)}\frac{1-|z|^2}{R}\int_{B((z,\sqrt{1-|z|^2}e^{i\theta}),CR)}w_ld\sigma=\frac1{R^{n+1}}\int_{B((z,\sqrt{1-|z|^2}e^{i\theta}),CR)}w_ld\sigma,
\end{split}\end{equation*}
with a simmilar relationship for $w^{-(p'-1)}$. Since $w_l$ is an $A_p$ weight in ${\bf S}^{n+1}$, we are done.\qed

The following result gives examples of $A_p(\B^n)$ weights obtained from weights in ${\bf S}^n$.
\begin{lemma}\label{exampleweights}
Assume $w\in A_p({\bf S}^{n})$. Then  the weight defined by
$${\widetilde{w}}(z)=\frac1{(1-|z|^2)^n}\int_{I_z}w(\zeta)d\sigma(\zeta),$$ $z\in \B^n$, where $I_z=\{\zeta\in{\bf S}^n\,;\,|1- \frac{z}{|z|}\overline{\zeta}|\leq c(1-|z|^2)\}$, $c>0$, is in $A_p(\B^n)$.
 \end{lemma}
 {\bf Proof of lemma \ref{exampleweights}:}\par

We want to show that there exists $C>0$ such that if   $a=(a_1,\cdots, a_n)$, and $U=U(a,R)=\{\eta\in\B^n\,;\, \rho(\eta,a)<R\}$, then
$$
\frac1{v(U)}\int_{U}{\widetilde{w}}(z)dv(z) \left(\frac1{v(U)}\int_{U}{\widetilde{w}}^{-(p'-1)}(z)dv(z)\right)^{\frac{-1}{p'-1}}\leq C.
$$ 

As in the previous lemma, assume first that $1-|a|\leq \frac1{\delta}R$, $\delta>0$ to be chosen. We can reduce this case to the one where $a=(a_1,0,\cdots,0)\in{\bf S}^n$. We then have that if $z\in U$ and $\zeta\in I_z$, then $\zeta \in U(a,CR)$, for some fixed constant $C>0$. Indeed, the fact that $z=(z_1,\cdots,z_n)\in U(a,R)$ gives (see lemma \ref{pseudodistance}) that $\displaystyle{\sum_{i=2}^n |z_i|^2\preceq R}$ and $\displaystyle{||z_1|^2-|a_1|^2|\preceq (R+R^\frac12 (1-|a_1|^2)^\frac12)\simeq R}$. In particular, we deduce that $1-|z|\preceq R$, and
$$|1-\zeta\overline{a}|=\rho(\zeta,a)\preceq \rho(\zeta, z)+\rho(z,a)\preceq R.$$
Next. Fubini's theorem gives that if $D_\alpha(\zeta)=\{z\in\B^n;\, |1-z\overline{\zeta}|<\alpha(1-|z|^2)\}$,
\begin{equation*}\begin{split}&
\frac1{v(U)}\int_{U}{\widetilde{w}}(z)dv(z)= \frac1{v(U)}\int_{U}\frac1{(1-|z|)^n} \int_{I_z} w(\zeta)d\sigma(\zeta) dv(z)\\
&\preceq\frac1{R^{n+1}}\int_{B(a,CR)} w(\zeta)d\sigma(\zeta) \int_{U(a,R)\cap D_\alpha(\zeta)} \frac{dv(z)}{(1-|z|)^n}\simeq\frac{1}{R^n}\int_{B(a,R)} w(\zeta)d\sigma(\zeta),
\end{split}\end{equation*}
 where $D_\alpha(\zeta)=\{z\in\B^n\,;\, |1-z\overline{\zeta}|<c(1-|z|^2)\}$, and we have used that if $\zeta\in I_z$, $|1-z\overline{\zeta}|\preceq (1-|z|)$. 
 An analogous argument  to the one we have used in last lemma applied to $\widetilde{w}^{\frac{-1}{p-1}}$ finishes this case.
 
 If $R\leq \delta (1-|a|)$, we have that for any $z\in U$, $\widetilde{w}(z)\simeq \widetilde{w}(a)$,  and consequentely, the $A_p$ condition in this case is obvious. 
 \qed

\begin{remark} \label{triebel} This type of weights $\widetilde{w}$ appear in a natural way if one identifies weighted holomorphic Besov spaces with   weighted holomorphic Triebel-Lizorkin spaces $HF_s^{pq}$, $p=q$. 
  We recall, that if $s\geq 0$,  $[s]^+$ denotes the integer part of $s$ plus $1$ and $1<p<+\infty$, the weighted holomorphic Triebel-Lizorkin space $HF_s^{pp}(w)$ is the space of holomorphic functions $f$ in $\B^n$ for which
\begin{equation*}\begin{split}&||f||_{HF_s^{pp}(w)}=\\&\left( \int_{{\bf S}^n}\left(\int_{D_\alpha(\zeta)} |(I+R)^kf(z)|^q( 1-|z|^2)^{(k-s)p-n-1}dv(z)\right) w(\zeta) d\sigma(\zeta)\right)^\frac{1}{p}<+\infty,
\end{split}\end{equation*}  
where $I$ denotes the identity operator. Fubini's theorem gives then that
$$||f||_{HF_s^{pp}(w)}\simeq \int_{\B^n}|(I+R)^kf(z)|^q( 1-|z|^2)^{(k-s)p-1}\widetilde{w}(z)dv(z).$$
 \end{remark}

\begin{definition} A weight $w$ in ${\bf S}^{n+1}$ is in $D_\tau({\bf S}^{n+1})$  if there exists $C>0$ such that for any $B(\zeta,R)=\{\eta\in{\bf S}^{n+1},\, |1-\eta\overline{\zeta}|<R\}$, $\zeta\in{\bf S}^{n+1}$, $R>0$, and $j\geq 0$, $w(B(\zeta,2^jR))\leq C2^{j\tau} w(B(\zeta,R))$.  
\end{definition}

Observe that any doubling weight $w$ in ${\bf S}^{n+1}$, i.e.  a weight $w$ there exists $k>0$ such that $w(B(\zeta,2R))\leq k w(B(\zeta,R))$, is in $D_\tau ({\bf S}^{n+1})$ for $\tau=\frac{\log k}{\log 2}$. So the fact that a weight is in $D_\tau({\bf S}^{n+1})$ is related to the size of the doubling constant $k$.

   In the proof of lemma 2.3 we have  seen in fact that if $w$ is a weight in $\B^n$ and $w_l$ is the corresponding lifted weight in ${\bf S}^{n+1}$, then $w_l(B(z,R))\simeq \frac{R}{R+(1-|z'|)}w_l(U(z',R))$, where $z=(z',z_{n+1})\in{\bf S}^{n+1}$. It is then natural to define a weight $w$ in $D_\tau(\B^n)$ as follows.

\begin{definition}
 We say that a weight $w$ in $\B^n$ is in $D_\tau(\B^n)$  for some $\tau$, if there exists  $C>0$, such that for any $k\geq 1$, $z\in\B^n$ and $R>0$, \begin{equation}\label{doublingtau}w(U(z,2^kR))\leq C \frac{(1-|z|^2)+2^kR}{(1-|z|^2)+R} 2^{k(\tau-1)} w(U(z,R)).\end{equation}
 
Since $v(U(z,R))\simeq R^{n}(R+(1-|z|^2))$, this  condition $D_\tau(\B^n)$ can be rewritten as
$$\frac{w(U(z,2^jR))}{v(U(z,2^jR))}\leq C \frac{2^{j\tau}w(U(z,R)) }{2^{j( n+1)}v(U(z,R))} .$$
\end{definition}
In many occassions, it  happens that in a natural way condition (\ref{doublingtau}) is only satisfied for those integers $k$ such that the ball $U(z,2^kR)$ touches the boundary of $\B^n$. We then have the following definition:
\begin{definition}
 We say that a weight $w$ in $\B^n$ is in $d_\tau(\B^n)$  for some $\tau$, if there exists  $C>0$, such that for any $k\geq 1$, $z\in\B^n$ and $R>0$ satisfying that $U(z,2^kR)$ touches ${\bf S}^n$, \begin{equation}\label{doublingtau1}w(U(z,2^kR))\leq C \frac{R}{(1-|z|^2)+R} 2^{k \tau } w(U(z,R)).\end{equation}
 \end{definition}

 As we have already observed, we have the following lemma.
 \begin{lemma}\label{justification}
 A weight $w$ is in $ D_\tau(\B^n)$, $\tau>0$, if and only if the lifted weight $w_l$ is in $D_\tau({\bf S}^{n+1})$.
 \end{lemma}
 {\bf Proof of lemma \ref{justification}:}\par
   Let $z=(z',z_{n+1})\in{\bf S}^{n+1}$, $R>0$ and $j\geq 1$. If $1-|z'|^2\leq R$,  $$w_l(B(z,2^jR))\simeq w(U(z',2^jR))\preceq 2^{j \tau}w(U(z',R))\simeq 2^{j\tau}w_l(B(z,R)).$$

If  $1-|z'|^2> 2^jR$, we have that (see the proof of lemma \ref{lifted}) 
$w_l(B(z,2^jR))\simeq \frac{2^jR}{1-|z'|^2}w(U(z', 2^jR))$, and $w_l(B(z,R))\simeq \frac{R}{1-|z'|^2}w (U(z',  R))$. Hence, 
\begin{equation*}\begin{split}&w_l(B(z,2^jR))\simeq \frac{2^jR}{1-|z'|^2}w(U(z', 2^jR))\\&
\preceq\frac{2^jR}{1-|z'|^2}2^{j\tau}w(U(z',  R))\simeq 2^{j\tau}w_l(B(z,R)).
\end{split}\end{equation*}

If $ 1-|z'|^2>R$, and $1-|z'|^2\leq 2^jR$, we have that (see Lemma \ref{lifted}) 
$w_l(B(z,2^jR))\simeq  w(U(z', 2^j R))$, and $w_l(B(z,R))\simeq \frac{R}{1-|z'|^2}w(U(z',  R ))$. Hence
$$w_l(B(z,2^jR))\simeq w(U(z', 2^j R))\preceq 
\frac{2^jR}{1-|z'|^2}2^{j\tau} w(U(z',R))\preceq 2^{j\tau}
  w_l(B(z,R)).
$$
 The other implication is proved in a simmilar way.\qed

Observe that if $w\equiv 1$, i.e., if $w$ is Lebesgue measure on $\B^n$, then $w\in D_\tau$, $\tau=n+1$. The following simple lemma will show that without loss of generality we always may assume that $\tau\geq n+1$.

\vspace{20pt}

\begin{lemma}\label{taugeqn+1}
If   $w$  is a non identically zero weight in $ L^1(\B^n)$ which is in $d_\tau$, then $\tau\geq n+1$.
\end{lemma} 
{\bf Proof of lemma \ref{taugeqn+1}:}\par
Assume $\tau <n+1$, and let $U\subset \B^n$ a ball in $\B^n$ that touches ${\bf S}^n$. The doubling condition on $w$ gives that for any $k\geq 1$, ${\displaystyle w(U)\preceq 2^{k\tau}w(2^{-k}U)}$. Consequently,
\begin{equation}\label{contradiction}
2^{k(n+1-\tau)}w(U)\preceq \frac{w(2^{-k}U)}{2^{-k (n+1) }}.\end{equation}
The differentiation theorem (see for instance Theorem 5.3.1 in \cite{rudin}) gives   that  for almost $z\in\B^n$, ${\displaystyle \lim_{k\rightarrow +\infty}\frac{w(2^{-k}U)}{2^{-k n }}\simeq 1.}$ Since we are assuming that $\tau<n+1$, this gives a contradiction with (\ref{contradiction}).\qed

 \begin{definition} \label{definitionBp}
We say that a weight $\omega$ is in $B_{p}(\B^n)$  (see \cite{bekolle}) if there exists $C>0$ such that for any ball $U(z,R)$  that touches ${\bf S}^n$, i.e., $\overline{U(z,R)}\cap {\bf S}^n\neq \emptyset$,  
$$\frac1{v(U(z,R))}\int_{U(z,R)}\omega dv \left( \frac1{v(U(z,R))}\int_{U(z,R)}\omega^{-(p'-1)}dv\right)^{p-1} \leq C.$$
\end{definition}
Obviously, any   $A_p(\B^n)$ weight satisfies the condition $B_p(\B^n)$.

We next observe that any weight $w\in B^p(\B^n)$ is in $d_\tau$ for $\tau=p(n+1)$. 

\begin{lemma}\label{doublingBp}
Let $1<p<+\infty$ and $w\in B^p(\B^n)$. There exists  $C>0$, such that for any $k\geq 1$, $z\in\B^n$ and $R>0$ satisfying that $U(z,2^kR)$ touches ${\bf S}^n$, $$w(U(z,2^kR))\leq C\left( \frac{R}{(1-|z|^2+R)}\right)^p 2^{k\tau}w(U(z,R)) \leq C\frac{R}{(1-|z|^2)+R} 2^{k \tau} w(U(z,R)),$$ 
 where $\tau=p(n+1)$.\end{lemma} 
{\bf Proof of lemma \ref{doublingBp}:}\par
Let $z\in\B^n$, $R>0$, $k\geq1$ such that $U(z,2^kR)$ touches ${\bf S}^n$. We then have that by H\"older's inequality,
\begin{align*}&
v(U(z,R))=\int_{U(z,2^kR)}\chi_{U(z,R)}(y)dv(y)=\int_{U(z,2^kR)}\chi_{U(z,R)}(y) w(y)^\frac1{p}w(y)^{-\frac1{p}}dv(y)\\& \leq\left(\int_{U(z,R)}\omega(y)dv(y)\right)^\frac1{p} \left(\int_{U(z,2^kR)}\omega^{-(p'-1)}(y)dv(y)\right)^\frac{p-1}{p}\\
&\preceq\left( \int_{U(z,R)}\omega(y)dv(y)\right)^\frac1{p} \left(\int_{U(z,2^kR)}\omega(y)dv(y)\right)^{-\frac1{p}} v(U(z,2^kR)),
\end{align*}
where in last inequality we have used the fact that $w\in B_p(\B^n)$. Thus we deduce that
\begin{equation}\label{doublinginequality}
w(U(z,2^kR))\leq C\left( \frac{v(U(z,2^kR))}{v(U(z,R))}   \right)^p w(U(z,R)).
\end{equation}
But lemma \ref{pseudodistance} together with the fact that $U(z,2^kR)$ touches ${\bf S}^n$, gives that
$$
\frac{v(U(z,2^kR))}{v(U(z,R))}\simeq \frac{(2^kR)^n(2^kR+(1-|z|^2)}{R^n(R+(1-|z|^2)}\preceq 
\frac{(2^kR)^{n+1}}{R^n(R+(1-|z|^2)}\preceq 2^{k(n+1)}\frac{R}{R+(1-|z|^2)}.
$$
Plugging the above inequality in (\ref{doublinginequality}), we deduce that
$$w(U(z,2^kR))\leq C \left(\frac{ R}{(1-|z|^2)+R} \right)^p 2^{k\tau} w(U(z,R)),$$
with $\tau=(n+1)p$.\qed
 
 As a  consequence of last lemma, we have an equivalent definition of $B_p(\B^n)$ weights which coincides with the weights in $B_p(\B^n)$ introduced in \cite{bekolle}:  a weight $w$ is in $B_p(\B^n)$ if there exists $C>0$ such that for any tent $T(\zeta,R)=\{z\in \B^n\,;\, |1-z\overline{\zeta}|\leq R\}$, $\zeta\in{\bf S}^n$,
 $$\frac1{v(T(\zeta,R))}\int_{T(\zeta,R)}\omega dv \left( \frac1{v(T(\zeta,R))}\int_{T(\zeta,R)}\omega^{-(p'-1)}dv\right)^{p-1} \leq C.$$
 This observation is a consequence of the fact that if a ball $U(\zeta,R)$ touches ${\bf S}^n$, then it is included in a tent of radious comparable to $R$, and conversely, a tent of radious $R$ is included in a ball that touches ${\bf S}^n$ of comparable radious.
 
 The weights in $B_p(\B^n)$ are characterized as the ones for which the Bergman  projector $B$ given by $$Bf(z)=\int_{\B^n} \frac{f(y)}{(1-z\overline{y})^{n+1}}dv(y),$$
 is a continuous operator from $L^p(w)$ to itself (see \cite{bekolle}).

 Let's  give some  examples of  weights in $B_p(\B^n)$ or in $ d_\tau(\B^n)$.
\begin{proposition}\label{exampledoublingweight}
Let $1<p<+\infty$, and let  $\varphi:(0,1]\rightarrow\R $ be a nonnegative monotone function, $C>0$, $\alpha>0$, satisfying one of the following alternative assumptions:
\begin{itemize}\item[(i)]     $\varphi$ is nondecreasing, and
$\displaystyle{ \varphi(2^kx)\leq C 2^{\alpha k }\varphi(x)}$.

\item[(ii)]    $\varphi$ is nonincreasing, and $\displaystyle{ \varphi( x)\leq C2^{k\alpha}\varphi(2^k x)}$ . 
\end{itemize} Let $w_\varphi(z)=\varphi(1-|z|)$. We then have:

\begin{itemize}
\item[(a)] The weight $w_\varphi$ is in $ B_p(\B^n)$ if and only if $0<\alpha<p-1$ if $\varphi$ is  in case (i) or $0<\alpha<1$ if it is in case (ii).
\item[(b)]  The weight $w_\varphi$ is in $ d_\tau$ if    $\varphi$ is case (i) and $\tau > n+\alpha+1$ or  $\varphi$ is in case (ii) and $\tau\geq n+1$.

\end{itemize}

 \end{proposition}
{\bf Proof of proposition \ref{exampledoublingweight}:}\par
Fubini's theorem gives that \begin{equation}\label{formula1}w_\varphi(U(z,R))=\int_{U(z,R)} \varphi(1-|y|) dv(y)\simeq R^{n-1} (R+R^\frac12(1-|z|)^\frac12) \int_{(1-|z|)}^{(1-|z|)+(R+R^\frac12 (1-|z|)^\frac12) }\varphi(t) dt.\end{equation}

We will show that if $\varphi$ satisfies the hypothesis in (a), then the weight $w_\varphi$ is in $B_p(\B^n)$. If $T(z,R)$ is a tent, (\ref{formula1}) gives that
$w_\varphi(T(z,R)) \simeq R^n \int_0^R \varphi(t)dt$, and consequently, it is enough to show that
$$\frac1{R}\int_0^R \varphi(t)dt \left( \frac1{R} \int_0^R\varphi^{-(p'-1)} (t)dt\right)^\frac1{p'-1} \leq C.$$

Assume first that $\varphi$ is nondecreasing.  
We   then have that for any $R>0$, $$\varphi(\frac{R}2)R\leq \int_{\frac{R}2}^R \varphi(t) dt\leq \int_0^R \varphi(t) dt\leq R\varphi(R).$$
Since  $\varphi(2x)\simeq \varphi(x)$, we obtain that ${\displaystyle{\int_0^R \varphi(t) dt\simeq R\varphi(R)}}$.

Next,
$$
\int_0^R \varphi^{-(p'-1)} (t)dt =\sum_{k=0}^{+\infty} \int_{\frac{R}{2^{k+1}}}^{\frac{R}{2^{k}}}\varphi^{-(p'-1)} (t)dt\simeq  
\sum_{k=0}^{+\infty} 2^{-k}R \varphi(\frac{R}{2^{k}})^{-(p'-1)}.
$$
Thus, it suffices to check that 
$$\sum_{k=0}^{+\infty} 2^{-k}  \varphi(\frac{R}{2^{k}})^{-(p'-1)}\preceq \frac1{\varphi(R)^{p'-1}}.$$
If we denote $a_{k+m}=2^{-k}\varphi(\frac{R}{2^{k}})^{-(p'-1)}$, the above is equivalent to show that
$$\sum_{k=0}^{+\infty} a_{m+k} \preceq a_m,$$ for any $m\geq 0$. But such condition it turns to be equivalent to
$$a_{m+k}\preceq \frac{a_m}{(1+\delta)^k},$$
for some $\delta>0$ (see for example subsection 5.4 in \cite{arcozzirochbergsawyer}). And that is a restatement of the doubling condition satisfied by $\varphi$.

If $\varphi$ is nonincreasing and such that $\varphi(2x)\simeq \varphi(x)$,  $\varphi^{-(p'-1)}$ is nondecreasing and satisfies that $\varphi^{-(p'-1)}(2x)\simeq \varphi^{-(p'-1)}(x)$. Hence, the above argument shows that if 
$$\varphi^{-(p'-1)}(2x)\preceq (2-\varepsilon)^{k(p'-1)}\varphi^{-(p'-1)}(x),$$
i.e., if $\varphi(x)\preceq (2-\varepsilon)^{k} \varphi(2x)$ then $w_{\varphi^{-(p'-1)}}(z)=\varphi(1-|z|)^{-(p'-1)}\in B_{p'}(\B^n)$. And that   is equivalent to say that $w\in B_p(\B^n)$.

The other implication is proved in a simmilar way.

We now prove (b). 
We want to show that under the conditions in (b), $w_\varphi\in d_\tau(\B^n)$, that is
$$w_\varphi (U(z,2^jR))\leq C \frac{(1-|z|^2)+2^jR}{1-|z|^2+R}2^{j(\tau-1)} w_\varphi (U(z,R)),$$ for any $U(z,R)$ and $j\geq 0$ such that $U(z,2^jR)$ touches ${\bf S}^n$.

We first recall  that 
   if $\displaystyle{M\leq \frac{x}2}$, we have that for any $t\in[x-M,x+M]$, $\frac{x}2\leq t\leq \frac{3x}2$, and consequently, $\displaystyle{ \int_{x-M}^{x+M} \varphi(t) dt\simeq \varphi(x) M}$. If on contrary, $M\geq \frac{x}2$, then $\int_{x-M}^{x+M}\varphi(t)dt\simeq M\varphi(M)$.
 
 Let $z\in \B^n$, $R>0$ and $j\geq 1$. We consider two different possibilities:
\begin{itemize}
 
\item[(1)] $\frac{(1-|z|^2)}2\leq R$.
\item[(2)] $R\leq \frac{(1-|z|^2)}2\leq 2^jR$.
\end{itemize}
  If $\frac{(1-|z|^2)}2\leq R$, (\ref{formula1}) and the above considerations give easily that
$w_\varphi(U(z,2^jR))\simeq(2^jR)^{n+1}\varphi(2^jR)$, and $w_\varphi(U(z,R)) \simeq R^{n+1}\varphi(R)$. Thus the  condition $d_\tau(\B^n)$ is fulfilled provided
$$(2^jR)^{(n+1)}\varphi(2^jR)\preceq 2^j2^{j(\tau-1)}R^{n+1}\varphi(R),$$   condition  that is equivalent to
$$\varphi(2^jR)\preceq 2^{j(\tau-n-1)}\varphi(R).$$
And the conditions on $\tau$ and $\varphi$ give that this estimate is satisfied .

If $R\leq \frac{(1-|z|^2)}2\leq 2^jR$. An argument analogous to   case (1), gives now that $w_\varphi(U(z,R))\simeq R ^{n }(1-|z|^2)\varphi(1-|z |)$. Hence it is enough to check in this case that
$$(2^jR)^{(n+1)}\varphi(2^jR)\preceq \frac{2^jR}{(1-|z|^2)}2^{j(\tau-1)} R^n (1-|z|^2)\varphi(1-|z |),$$
i.e., $  \varphi(2^jR)\preceq 2^{j(\tau-n-1)} \varphi(1-|z |)$. And this estimate is a consequence on the hypothesis on $\varphi$. \qed

 \begin{corollary}\label{examplealpha}
 If $w_\alpha(z)=(1-|z|)^\alpha$, $-1<\alpha<p-1$, then the weight $w_\alpha$ is in $ B_p(\B^n)$. If   $0\leq \alpha<p-1$ and $\tau = n+\alpha+1$, or if $-1<p<0$ and $\tau=n+1$, then $w_\alpha\in d_\tau(\B^n)$.
 \end{corollary}
 
 The techniques we will apply in order to work with weighted holomorphic Besov spaces in $\B^n$ require that the weights $w$  are in $A_p(\B^n)\cap D_\tau(\B^n)$. The purpose of the following pair of technical results is to show that we can in fact weaken these conditions and impose that the weight $w$ is in the bigger class $B_p(\B^n)\cap d_\tau(\B^n)$. The way to achieve this is via the regularisations of the weights.
 
\begin{definition}\label{definitionregularisation}
If $w$ is a weight in $\B^n$, $0< \varepsilon<1$ and $U_\varepsilon (z)=U(z,\varepsilon (1-|z|^2))$, we define
$$
R_\varepsilon{w}(z)=\frac1{v(U_\varepsilon (z))}\int_{U_\varepsilon(z)} w(\eta)dv(\eta).
$$
\end{definition}
\begin{remark}
As an immediate consequence of lemma \ref{pseudodistance} and lemma \ref{doublingBp}, we have that all the regularisations are equivalent, that is, if $\varepsilon,\varepsilon'>0$, $R_\varepsilon w(z)\simeq R_{\varepsilon'}w(z)$, for any $z\in\B^n$, with constants that do not depend on $z$. We just have to observe that if $\varepsilon>0$ is fixed, there exists $C>\varepsilon$ such that $U(z,C(1-|z|)$ touches ${\bf S}^n$. The fact that $w$ satisfies a doubling condition gives that $w(U_\varepsilon)\simeq w(U(z,C(1-|z|^2)))$.

 Observe that the regularisation of a weight $w$ satisfies that $R_\varepsilon(R_\varepsilon w)\simeq R_\varepsilon w$.
\end{remark}

  It is worthwhile to recall that analogous regularisations where already considered among others by \cite{bekolle} and \cite{luecking},  where the balls $U_\varepsilon^d(z)=\{\eta\in\B^n\,;\, d(z,\zeta)<\varepsilon(1-|z|^2)\}$ were defined with respect to the pseudodistance $d(z,\zeta)=||z|-|\zeta||+|1-(z\overline{\zeta})/|z||\zeta||$.

\begin{lemma}\label{regularisation}
Let $1<p<+\infty$ and assume that $w$ is a  weight in $B_p(\B^n)$. Then the weight
$$
R_\varepsilon w(z)=\frac1{v(U_\varepsilon(z))}\int_{U_\varepsilon(z)} w(\eta)dv(\eta),
$$
is in $A_p(\B^n)$.
\end{lemma}
{\bf Proof of lemma \ref{regularisation}:}\par

We want to show that there exists $C>0$
 such that for any  ball $U=U(a,R)=\{\eta\in\B^n\,;\, \rho(\eta,a)<R\}$ associated to the pseudodistance $\rho$, 
$$\left( \frac{1}{v(U)}\int_U R_\varepsilon w dv\right)\left( \frac{1}{v(U)}\int_U (R_\varepsilon w)^{-(p'-1)}dv \right)^{\frac1{p'-1}}\leq C.$$

As we have already observed,  without loss of generality  we may assume that $\varepsilon>0$ is small enough, since for every $\varepsilon,\varepsilon'>0$, then $R_\varepsilon w\simeq R_{\varepsilon'}w$.

Suppose first that $\delta(1-|a|^2)\leq R$, $\delta>0$ to be chosen later on. In this case, Lemma \ref{pseudodistance} gives that $v(U)\simeq R^{n+1}$. Since we also have that in that case $U(a,R)$ is included in a ball in ${\overline \B}^n$ centered at a point in ${\bf S}^n$ of radious comparable to $R$, we also may assume without loss of generality that $a\in {\bf S}^n$, and that $U=U(a,CR)$, $C>0$. In particular we have that for any $\eta\in U(a,R)$, $1-|\eta|^2\preceq R$.
Thus if $z\in U(a,CR)$ and $y\in U_\varepsilon(z)$, $y\in U(a,CR)$, and Fubini's theorem gives that
\begin{equation} \label{Ap}\int_{U(a,CR)} R_\varepsilon w(z)dv(z)\simeq \int_{U(a,CR)}\frac1{(1-|z|^2)^{n+1}}\int_{U_\varepsilon(z)}w(y)dv(y)dv(z).\end{equation}
Next, if $y\in U_\varepsilon(z)$, and $\varepsilon >0$ is small enough,  there exists $\varepsilon'<1$ such that $z\in U_{\varepsilon'}(y)$, and $1-|z|^2\simeq 1-|y|^2$. Thus (\ref{Ap}) is bounded by
$$\int_{U(a,CR)}w(y)\frac1{(1-|y|^2)^{n+1}} \int_{U_{\varepsilon'}(y)} dv(z) dv(y)\preceq \int_{U(a,CR)} w(y)dv(y).
$$

In order to estimate the integral involving $(R_\varepsilon w)^{-(p'-1)}$, we use the fact that $w\in B^p(\B^n)$   and H\"older's inequality to get that,
$$\frac{1}{v(U_\varepsilon(z))}\int_{U_\varepsilon (z)}w(y)dv(y) \left(\frac{1}{v(U_\varepsilon(z))}\int_{U_\varepsilon(z)} w^{-(p'-1)}  \right)^\frac1{p'-1}\simeq 1.$$
Consequently   $$R_\varepsilon w\simeq \left( R_\varepsilon(w^{-(p'-1)})\right)^\frac{-1}{p'-1}.$$
This gives that
\begin{align*}&\frac1{v(U(a,CR))}\int_{U(a,CR)} \left( R_\varepsilon w(z)\right)^{-(p'-1)}dv(z) \\&\simeq 
\frac1{v(U(a,CR))} \int_{U(a,CR)} R_\varepsilon(w^{-(p'-1)})(z)dv(z),\end{align*}
and the argument in (\ref{Ap}) applied to $R_\varepsilon(w^{-(p'-1)})$  together with the fact that $w\in B_p(\B^n)$, gives that in case $1-|a|\leq \frac1{\delta}R$, then 
$$\frac1{v(U(a,R))}\int_{U(a,R)}R_\varepsilon w(z)dv(z)\left( \frac1{v(B(a,R))}\int_{B(a,R)} (R_\varepsilon w)^{-(p'-1)}(z)dv(z) \right)^\frac1{p'-1}\leq C.$$
Assume next that $R\leq \delta(1-|a|^2)$, and $\delta$ is small enough. In that case, for any $z\in U(a,R)$, $1-|z|^2\simeq 1-|a|$, and consequently $R_\varepsilon w(z)\simeq R_\varepsilon w(a)$  for any $z  \in U(a,R)$. This is a direct consequence of the fact that $w$ is doubling and that if $z\in U(a,R)$, there exists $\varepsilon',\varepsilon''>0$ such that
$U_\varepsilon(z)\subset U_{\varepsilon '}(a)\subset U_{\varepsilon''}(z)$.
\qed

 \begin{lemma}\label{regdobl}
 If $w$ is a doubling weight in $\B^n$, its regularisation $R_\varepsilon w$ also satisfies a doubling condition.
 \end{lemma}
{\bf Proof of lemma \ref{regdobl}:}\par

 Assume that there exists $C>0$ such that for any $z\in \B^n$, and $R>0$,
 $$w(U(z,2R))\leq Cw(U(z,R)).$$ 
 We want to check that $R_\varepsilon w$ satisfies a simmilar condition. Let $z\in\B^n$ and $R>0$, and assume first that $R  \leq \delta(1-|z|^2)$, with $\delta>0$ small enough so that for any $y\in U(z,2R)$, $1-|y|^2\simeq 1-|z|^2$. 
 We then have that
 \begin{equation*}\begin{split}&R_\varepsilon w(U(z,  2  R))= \int_{U(z,  2  R)}R_\varepsilon w(\eta)dv(\eta)\simeq  R_\varepsilon w(z) v(U(z,  2 R ))\\&\simeq
R_\varepsilon w(z) R^n (1-|z|^2)\simeq 
  \int_{U(z, R)}R_\varepsilon w(\eta)dv(\eta) .\end{split}\end{equation*}

If $\delta(1-|z|)\leq R$, Fubini's theorem and the fact that $w$ satisfies a doubling condition  give
\begin{equation*}\begin{split}&R_\varepsilon w(U(z,  2  R))= \int_{U(z,  2  R)} \frac1{(1-|z|^2)^{n+1}} \int_{U_\varepsilon(z)} w(y)dv(y) dv(z)   \\&\preceq
\int_{U(z,  C  R)} w(y) \frac{1}{(1-|z|^2)^{n+1}}\int_{U_{\varepsilon'}(y)}dv(z)dv(y)\preceq \int_{U(z,    R)} w(y) dv(y) \\&\simeq
\int_{U(z,    R)} w(y) \frac{1}{(1-|z|^2)^{n+1}}\int_{U_{\varepsilon'}(y)}dv(z)dv(y)\simeq R_\varepsilon w(U(z,   R)).\qed
\end{split}\end{equation*}

\begin{proposition}\label{doublingregularisation}
Let $1<p<+\infty$ $\tau\geq n+1$, and assume that $w$ is a   weight   satisfying that $w\in d_\tau(\B^n)$. Then the weight
$ 
R_\varepsilon w $
is in $D_\tau(\B^n)$.
\end{proposition}
{\bf Proof of proposition \ref{doublingregularisation}:}\par

The hypothesis on $w$ gives that for any $z\in \B^n$, $j\geq 0$ and $R>0$, such that $U(z,2^jR)$ touches ${\bf S}^n$, then
$$
w(U(z,2^jR))\preceq \frac{(1-|z|^2)+2^jR}{(1-|z|^2)+R}2^{j(\tau-1)}w(U(z,R)).$$
In order to check that $R_\varepsilon w\in D_\tau(\B^n)$,  given $z\in\B^n$, $j\geq 0$ and $R>0$, we will consider the following three possibilities:
\begin{itemize}
\item[(a)] $\displaystyle{    2^{j-1} R\leq \delta (1-|z|^2)}$.
\item[(b)]$\displaystyle{   R< \delta (1-|z|^2)\leq  2^{j-1} {R}}$.
\item[(c)] $\displaystyle{     \delta  (1-|z|^2)< R}$.
\end{itemize}
Here $\delta>0$ is some fixed constant to be chosen later on. 

We begin with case (a). Our first observation is that in that case $R_\varepsilon w(y)\simeq R_\varepsilon(z)$ for $y\in U(z,   2^{j-1} R)$. This is an immediate consequence of the doubling condition on $w$ and the fact that if $y\in U(z,  2^{j-1} R)$, $\rho(y,z)<< (1-|z|^2)$, and hence $1-|y|^2\simeq 1-|z|^2$.

Hence, by the preceding lemma,   
\begin{equation*}\begin{split}&R_\varepsilon w(U(z,  2^{j} R))\simeq R_\varepsilon w(U(z,  2^{j-1} R))= \int_{U(z,  2^{j-1} R)}R_\varepsilon w(\eta)dv(\eta)\simeq  R_\varepsilon w(z) v(U(z,  2^{j-1}R ))\\&\simeq
R_\varepsilon w(z) (2^jR)^n (1-|z|)\simeq 
 2^{jn}\int_{U(z, R)}R_\varepsilon w(\eta)dv(\eta)\preceq 2^{j(\tau-1)}R_\varepsilon w(U(z, R)) ,\end{split}\end{equation*}
where in last inequality we have used that  $\tau\geq n+1$.This shows case (a).
We consider now case (b). Let $j_0\geq 1$ such that $2^{j_0-1}R\leq \delta (1-|z|^2)<2^{j_0}R$.

Fubini's theorem, and the fact that $w\in d_\tau(\B^n)$ gives that
\begin{equation}\label{estimate3}\begin{split}& \int_{U(z,2^jR)} R_\varepsilon w(y)dv(y)\simeq \int_{U(z,C2^{j}R)}\frac1{(1-|y|^2)^{n+1}}\int_{U(y,\delta(1-|y|^2))} w(\eta)dv(\eta) dv(y)\\&\preceq
\int_{U(z,C2^jR)}w(\eta)dv(\eta)\preceq \frac{2^{j-j_0}R+(1-|z|^2)}{2^{j_0}R+(1-|z|^2)}2^{(j-j_0)(\tau-1)} \int_{U(z,2^{j_0}R)}w(\eta)dv(\eta).
\end{split}\end{equation}
Since $R\leq \delta (1-|z|^2)<2^{j_0}R$, the argument stablished in case (a) gives that $R_\varepsilon w$ is ''frozen'' on $U(z,\frac{R}2)$. This observation, together with the fact that $R_\varepsilon w$ satisfies a doubling condition, gives that
\begin{equation*}\begin{split}& 
\int_{U(z,R)}R_\varepsilon w(\eta)dv(\eta)\simeq \int_{U(z,\frac{R}2)}R_\varepsilon w(\eta)dv(\eta)\simeq R_\varepsilon w(z) v(U(z,\frac{R}2))\\
&\simeq R_\varepsilon w(z) R^n(\frac{R}2+(1-|z|^2))\simeq \frac{R^n 2^{j_0}R}{(1-|z|^2)^{n+1}}\int_{U(z,2^{j_0}R)}w(\eta)d(\eta),
\end{split}\end{equation*}
where in last estimate we have used that $R\leq \delta(1-|z|^2)\simeq 2^{j_0}R$.
Consequently, if we plug the above calculation in (\ref{estimate3}), we deduce that in order to prove the doubling condition in case (b) it is enough to show that
$$\frac{2^{j-j_0}R+(1-|z|^2)}{2^{j_0}R+(1-|z|^2)}2^{(j-j_0)(\tau-1)}\frac{(1-|z|^2)^{n+1}}{R^n 2^{j_0}R}\preceq \frac{2^jR+(1-|z|^2)}{R+(1-|z|^2)}2^{j(\tau-1)}.$$
Now, the fact that we are in case (b) gives that $2^jR+(1-|z|^2)\simeq 2^jR$, and $R+(1-|z|^2)\simeq (1-|z|^2)\simeq 2^{j_0}R$, and consequently the above estimate can be rewritten equivalently as
$$\frac{2^{j-j_0}R+(1-|z|^2)}{2^{j_0}R}2^{-j_0(\tau-1)} \frac{(2^{j_0}R)^{n+1}}{2^{j_0}R^{n+1}} \frac{2^{j_0}R}{2^jR}=\frac{2^{j-j_0}R+(1-|z|^2)}{R}2^{j_0(n+1-\tau)}2^{-j} \preceq 1.$$
But the fact that $2^{j-j_0}R+(1-|z|^2)\leq 2^jR+(1-|z|^2)\simeq 2^jR$ gives that the left hand side of the above is bounded from above by
$$C\frac{2^jR}{R}2^{j_0(n+1-\tau)}2^{-j}=C2^{n+1-\tau}\preceq 1,$$
since $\tau\geq n+1$.

We finally have to deal with case (c), i.e. the case where $\delta(1-|z|^2)\leq R$. We   have that if $y\in U(z, 2^{j-1}R)$, and $\eta\in U(y,\delta(1-|z|^2))$, then $\eta\in U(z,C2^jR)$, and consequently, Fubini's theorem gives that
\begin{align}\label{uperestimate}& \int_{U(z,2^jR)} R_\varepsilon w(y)dv(y)\simeq \int_{U(z,2^{j-1}R)}\frac1{(1-|y|^2)^{n+1}}\int_{U(y,\delta(1-|y|^2))} w(\eta)dv(\eta) dv(y)\\&\preceq
\int_{U(z,C2^jR)}w(\eta)dv(\eta)\preceq 2^{j\tau}\int_{U(z,R)}w(\eta)dv(\eta),
\end{align}
where we have used that $w\in d_\tau$.

On the other hand, $R_\varepsilon w$  satisfies a doubling condition. Thus, if $M>0$ is fixed, Fubini's theorem gives that
\begin{align}\label{lowerestimate}&
\int_{U(z, R)}w(\eta)dv(\eta)\simeq \int_{U(z,\delta_1 R)}w(\eta)\frac1{(1-|\eta|^2)^{n+1}}\int_{U(\eta,\delta(1-|\eta|^2))}dv(y)dv(\eta) \\&\preceq\int_{U(z,MR)}\frac1{(1-|y|^2)^{n+1}}\int_{U(y,\delta(1-|y|^2))} w(\eta)dv(\eta) dv(y)\\&\simeq
\int_{U(z,MR)} R_\delta w(y)dv(y)\simeq \int_{U(z,R)} R_\varepsilon w(y)dv(y),\end{align}
 where in the second estimate we have used that since $\delta(1-|z|^2) \leq R$, then $(1-|\eta|^2)\preceq R$ for any $\eta\in U(z,R)$, and that if $M>0$ is big enough,  then for any $\eta\in U(z,R)$, $U(\eta, \delta(1-|\eta|^2))\subset U(z,MR)$. \qed
 
 \begin{remark}
 We have shown that  the regularisation $R_\varepsilon w$ of a weight $w$ in $\B_p$ is in the smaller class $A_p(\B^n)$. In particular,  $R_\varepsilon w$ satisfies a doubling condition and consequently, if $w\in B_p(\B^n)\cap d_\tau(\B^n)$, then $R_\varepsilon w\in A_p(\B^n)\cap D_\tau(\B^n)$.
 \end{remark}

We will show next that the weights $\widetilde{w}$  introduced in lemma \ref{exampleweights} that were obtained from weights in $A_p({\bf S}^n)$ are  in  $A_p(\B^n)\cap d_{\tau+1}(\B^n)$ if $w\in D_\tau({\bf S}^n)$.
\begin{lemma}\label{exampleweightsinsn}
Assume $w\in A_p({\bf S}^{n})\cap D_{\tau}({\bf S}^n)$, $\tau\geq n$. Then  the weight defined by
$${\widetilde{w}}(z)=\frac1{(1-|z|^2)^n}\int_{I_z}w(\zeta)d\sigma(\zeta),$$ $z\in \B^n$, is in $A_p(\B^n)\cap d_{\tau+1}$.
 \end{lemma}
 {\bf Proof of lemma \ref{exampleweightsinsn}:}\par
 By lemma \ref{exampleweights} we know that ${\widetilde{w}}\in A_p(\B^n)$, so we are left to show that ${\widetilde{w}}\in d_{\tau+1}(\B^n)$. 
 Let $U(a, 2^kR)$ be a ball in $\B^n$ that touches ${\bf S}^n$. We want to show that
 $$\widetilde{w}(U(a,2^kR))\preceq \frac{2^kR}{(1-|z|^2)+R}2^{k \tau} \widetilde{w}(U(a,R)).$$

  Fubini's theorem gives that
 \begin{equation*}\begin{split}& 
 \int_{U(a,2^kR)}{\widetilde{w}}(z)dv(z)= \int_{U(a,2^kR)}\frac1{(1-|z|^2)^n}\int_{I_z}w(\zeta)d\sigma(\zeta) dv(z)\\&\simeq \int_{B(\frac{a}{|a|},C2^kR)} w(\zeta) d\sigma(\zeta) \int_{U(a,2^kR)\cap D_\alpha(\zeta)} \frac{dv(z)}{(1-|z|^2)^n} \simeq 2^kR \int_{B(\frac{a}{|a|},2^kR)}w(\zeta)d\sigma(\zeta).
\end{split}\end{equation*}  
Assume first that $R> \delta (1-|a|^2)$, $\delta>0$ small enough. Then last argument  can be applied to $U(a,R)$, and we get
 $$\int_{U(a,R)}{\widetilde{w}}(z)dv(z)\simeq R \int_{B(\frac{a}{|a|},R)}w(\zeta)d\sigma(\zeta).$$
 Thus in that case, the doubling condition reduces to check that
 $$2^k R w(B(\frac{a}{|a|},2^kR))\preceq 2^{k\tau} Rw(B(\frac{a}{|a|},R)),$$
 which follows from the fact that $w\in D_\tau$.
 If $R<\delta(1-|a|^2)$ and $\delta$ is small enough, we have observed in previous lemmas that $\widetilde{w}$ is "frozen" in $U(a,R)$. Consequently, 
 $$\widetilde{w}(U(a,R)) \simeq v(U(a,R)) \widetilde{w}(a)\simeq R^n(1-|a|^2)\frac1{(1-|a|^2)^n}\int_{I_a} w(\zeta)d\sigma(\zeta).$$
This observation, together with the fact that $w$ is in $D_\tau$ gives then that
\begin{equation*}\begin{split}& 
\widetilde{w}(U(a,2^kR))\simeq 2^kR w(B(\frac{a}{|a|},2^kR))\preceq  2^kR  \left(\frac{2^kR}{(1-|a|^2)} \right)^\tau w(B(\frac{a}{|a|}, (1-|a|^2))) \\&\simeq
2^kR  \left(\frac{2^kR}{(1-|a|^2)} \right)^\tau \frac{(1-|a|^2)^n}{R^n(1-|a|^2) } \widetilde{w}(U(a,R))=\frac{2^kR}{(1-|a|^2)}2^{k\tau} \left(\frac{R}{(1-|a|^2)}\right)^{\tau-n}\widetilde{w}(U(a,R)).
\end{split}\end{equation*} 
Since in this case, $\frac{R}{(1-|a|^2)}\preceq 1$, we are done.\qed

\section{Weighted holomorphic Besov spaces}

We now introduce the weighted holomorphic Besov  spaces. Let $w$ be an $B_p$-weight in $\B^{n}$, $1<p<+\infty$, and $s\in \R$ . The space $B_{s,k}^{p}(w,\B^n)$ is the space of holomorphic functions in $\B^n$ for which
$$||f||_{B_{s,k}^{p}(w,\B^{n})}^p=\int_{\B^{n}} |(I+R)^kf(y)|^p (1-|y|^2)^{(k-s)p-1}w(y)dv(y)<+\infty,$$
where $k\in\N$, $k>s$.
In fact, the definition of the weighted holomorphic Besov spaces does not depend on $k>s$. This is the object of the following result.

\begin{theorem}\label{independenceparameter}
Let $1<p<+\infty$, $s\in \R$, $k_1>k_2>s$, and $w$ a  $B_p$-weight in $\B^n$. We then have that the following are equivalent:

(i) $f\in B_{s,k_1}^{p }(w,\B^n)$.

(ii) $f\in B_{s,k_2}^{p}(w,\B^n)$.

\end{theorem}
{\bf Proof of theorem \ref{independenceparameter}:}\par

Assume that (i) holds. The fact that $f\in B_{s,k_1}^{p }(w,\B^n)$ means that $(I+R)^{k_1}f(y) (1-|y|^2)^{k_1-s-\frac1{p}}$ is in $L^p(wdv)$. In particular, there exists $p_1<p$ such that $(I+R)^{k_1}f(y) (1-|y|^2)^{k_1-s-\frac1{p}}$ is in $L^{p_1}(dv)$ (see \cite{cascanteortega3}, Lemma 1.1). Consequentely,  the kernel $c_N(1-|z|^2)^N/(1-\overline{z}y)^{n+1+N}$, for $N$ big enough, and for an adequate constant $c_N>0$, is a reproducing kernel for the function $f$ and its derivatives. We then have  
$$
(I+R)^{k_2}f(y)= c_N\int_{\B^n} (I+R)^{k_1}f(z)(I+R_y)^{k_2-k_1}\frac{(1-|z|^2)^N}{(1- \overline{z}y)^{n+1+N}}dv(z),
$$
where the operator $(I+R)^{-m}$ has the following integral representation (see for instance \cite{ortegafabrega})
$$
(I+R)^{-m}f(z)=\frac1{\Gamma(m)}\int_0^1 \left( \log \frac1{r} \right)^{m-1} f(rz)dr.$$
Thus
 $$|(I+R)^{k_2}f(y)|\preceq \int_{\B^n} \frac{|(I+R)^{k_1}f(z)|(1-|z|^2)^N}{|1-\overline{z}y|^{n+1+N+k_2-k_1}}dv(z),$$
  and we have,
 
\begin{equation}\begin{split}\label{splitting}
&||f||_{B_{s,k_2}^{p }(w,\B^n)}=\\&\sup_{ ||\psi||_{L^{p'}(wdv)}\leq 1} |\int_{\B^n} (I+R)^{k_2}f(y)(1-|y|^2)^{k_2-s-\frac1{p}} \psi(y)w(y)dv(y)|\\&\preceq
\sup_{ ||\psi||_{L^{p'}(wdV)}\leq 1}\int_{\B^n} \int_{\B^n} |(I+R)^{k_1}f(z)|\frac{(1-|z|^2)^N(1-|y|^2)^{ k_2-s-\frac1{p}}}{|1-\overline{z}y|^{n+1+N+k_2-k_1}}\psi(y) w(y) dv(y)dv(z).
 \end{split}\end{equation}

We now check  that the mapping 
\begin{equation}\label{operadoracotat}T_{N,k_1,k_2}(\psi)= \int_{\B^n} \psi(y)\frac{(1-|z|^2)^{N- (k_1-s-\frac1{p})}(1-|y|^2)^{k_2-s-\frac1{p}}}{|1-\overline{z} y|^{n+1+N+k_2-k_1}}w(y)dv(y),
\end{equation}
is bounded from $L^{p'}(wdv)$ to $L^{p'}(w^{-(p'-1)})$.
Indeed, if we denote $\alpha(y)=\psi(y)w(y)$, this holds if and only if the mapping
$$
\alpha\rightarrow \int_{\B^n} \frac{(1-|z|^2)^{N- (k_1-s-\frac1{p})}(1-|y|^2)^{k_2-s-\frac1{p}}}{|1-\overline{z} y|^{n+1+N+k_2-k_1}}\alpha(y)dv(y),$$
is bounded from $L^{p'}(w^{-(p'-1)})$ to itself. 

Since  $${\displaystyle\frac{(1-|z|^2)^{N- (k_1-s-\frac1{p})}(1-|y|^2)^{k_2-s-\frac1{p}}}{|1-\overline{z} y|^{n+1+N+k_2-k_1}}\preceq \frac{(1-|y|^2)^{k_2-s-\frac1{p}}}{|1-\overline{z} y|^{n+1+k_2-s-\frac1{p}}},}$$ and $k_2-s-\frac1{p}>-1$, this is a consequence of proposition 2 in \cite{bekolle}, where it is shown that the operator $$T_{k_2-s-\frac1{p}}f(z)=\int_{\B^n}\frac{(1-|y|^2)^{k_2-s-\frac1{p}}f(y)}{|1-\overline{z} y|^{n+1+k_2-s-\frac1{p}}}dv(y)$$ sends $L^{p'}(w^{-(p'-1)})$ to itself, provided   the weight $w^{-(p'-1)}$ is in the class $B_{p'}$. 

We then have from (\ref{splitting}) that 
\begin{equation*}\begin{split}
 &||f||_{B_{s,k_2}^{p }(w,\B^n)}\preceq\\&
\sup_{ ||\psi||_{L^{p'}(wdv)}\leq 1}\int_{\B^n}   |(I+R)^{k_1}f(z)| (1-|z|^2)^{k_1-s-\frac1{p}}w^\frac1{p}(z)T_{N,k_1,k_2}(\psi)(z)w^{-\frac1{p}}(z)dv(z)\\&
\preceq ||f||_{B_{s,k_1}^{p}(w,\B^n)}\sup_{ ||\psi||_{L^{p'}(wdv)}\leq 1}||T_{N,k_1,k_2}(\psi)||_{L^{p'}(wdv}  \\&\preceq ||f||_{B_{s,k_1}^{p }(w,\B^n)}\sup_{ ||\psi w||_{L^{p'}(w^{-(p'-1)}dv)}\leq 1} ||\psi||_{L^{p'}(wdv)}=||f||_{B_{s,k_1}^{p}(w,\B^n)}.
\end{split}\end{equation*}

The other implication is proved in a simmilar way.\qed

 By Theorem \ref{independenceparameter}, the spaces $B_{s,k}^{p}(w,\B^n)$ do not depend on $k>s$, and from now on we will denote them simply by $B_s^{p}(w,\B^n)$. Observe that if $w\in A_p({\bf S}^n)$, and ${\displaystyle\widetilde{w}(z)=\frac1{(1-|z|^2)^n}\int_{I_z}w(\zeta)d\sigma(\zeta)}$, we have that 
$B_{s }^{p}(\widetilde{w},\B^{n})=HF_s^{pp}(w)$.

Our next result shows that the weighted Besov space  associated to a weight in $B_s^p(\B^n)$ coincides with the corresponding weighted space of its regularisation. In particular, we deduce that in the definition of $B_s^p(w,\B^n)$ we can assume, without loss of generality that $w\in A_p(\B^n)$.

\begin{proposition}\label{weightedBesovregularisation}
Let $1<p<+\infty$, $s\in\R$, and assume that $w$ is a  weight in $B_p(\B^n)$. Then    
the spaces $B_s^p(w,\B^n)$ and $B_s^p(R_\varepsilon w,\B^n)$ coincide.
\end{proposition}
{\bf Proof of proposition \ref{weightedBesovregularisation}:}\par

Assume that $\varepsilon<1$ and $z\in\B^n$, an let $k>s$. The fact that $(I+R)^kf$ is holomorphic in $\B^n$ and that  $U_\varepsilon(z)$ is contained and contains an ellipsoid $E(z)$ in $\B^n$ of the same size,  gives immediately that
$$
|(I+R)^kf(y)|\preceq \frac{1}{|U_\varepsilon(y)|}\int_{U_\varepsilon(y)}|(I+R)^kf(z)|dv(z).
$$
On the other hand, $(1-|y|^2)\simeq (1-|z|^2)$ for any $y\in U_\varepsilon(z)$. Hence
$$(1-|y|^2)^\frac{(k-s)p-1}{p}|(I+R)^kf(y)|\preceq \frac{1}{|U_\varepsilon(y)|}\int_{U_\varepsilon(y)}(1-|z|^2)^\frac{(k-s)p-1}{p}|(I+R)^kf(z)|dv(z).
$$
Consequentely,
\begin{equation*}\begin{split} 
&||f||_{B_s^p(R_\varepsilon w)}=\int_{\B^n} |(I+R)^kf(y)|^p(1-|y|^2)^{(k-s)p-1} R_\varepsilon w(y)dv(y)\\
&\preceq\int_{\B^n}\left[\frac{1}{(1-|y|^2)^{n+1}}\int_{U_\varepsilon(y)}(1-|z|^2)^\frac{(k-s)p-1}{p} |(I+R)^kf(z)dv(z)|\right]^pR_\varepsilon  w(y)dv(y) \\
&\simeq\int_{\B^n}\left(R_\varepsilon\left( (1-|z|^2)^\frac{(k-s)p-1}{p}|(I+R)^kf(z)|  \right)(y)\right)^p
R_\varepsilon  w(y)dv(y).  \end{split}\end{equation*}

Since $w\in B_p(\B^n)$, 
$$
\int_{\B^n}(R_\varepsilon g)^pR_\varepsilon wdv\preceq \int_{\B^n} g^pwdv
$$
for any $g\geq 0$ (the proof is analogous to lemma 9 in \cite{bekolle}), we deduce that the above is bounded by
$$\int_{\B^n}\left((1-|y|^2)^{\frac{(k-s)p-1}{p}}|(I+R)^kf(y)|\right)^pw(y)dv(y)=C||f||_{B_s^p(w)}^p.$$
On the other hand, the fact that $|(I+R)^kf|^p$ is plurisubharmonic gives that
$$
|(I+R)^kf(y)|^p\preceq \frac{1}{|U_\varepsilon(y)|}\int_{U_\varepsilon(y)}|(I+R)^kf(z)|^pdv(z).
$$
Since in $U_\varepsilon(z)$, $(1-|y|^2)\simeq (1-|z|^2)$, we  have
\begin{equation*}\begin{split} 
&||f||_{B_s^p(w)}=\int_{\B^n} |(I+R)^kf(y)|^p(1-|y|^2)^{(k-s)p-1}   w(y)dv(y)\\
&\preceq\int_{\B^n} \left(R_\varepsilon\left(  (1-|y|^2)^\frac{(k-s)p-1}{p}(I+R)^kf      \right)^p(y)\right)w(y)dv(y)\preceq
\end{split}\end{equation*}

 Fubini's theorem gives that there exists $\varepsilon'>0$ such that for any $f,g\geq 0$, 
$$\int_\B^n f R_\varepsilon g dv \preceq \int_\B^n  g R_{\varepsilon'} f dv.
$$
Hence, the above is bounded by
$$\int_{\B^n}     (1-|y|^2)^ {(k-s)p-1}|(I+R)^k f (y)|^p  R_{\varepsilon'} w(y)dv(y)=C||f||_{B_s^p(R_\varepsilon w)}.\qed $$

  Our next goal is to study the relations between the weighted Besov spaces in $\B^n$ and the weighted Hardy-Sobolev spaces with respect to the lifted weight in ${\bf S}^{n+1}$. We will first show that the restriction operator maps $H^p_s(w_l)$ onto $B_{s-\frac{1}{p}}(w,\B^n)$.

We begin with a weighted restriction theorem. In order to make clearer the notation, in what follows we will write  $HF_s^{p,q}(w_l,\B^{n+1})$ instead of $HF_s^{p,q}(w_l)$. We recall (see \cite{cascanteortega3}) that $H_s^{p,2}(w_l,\B^{n+1})=H_s^2(w_l)$ and that if $q_0\leq q_1\leq +\infty$, ${\displaystyle HF_s^{pq_0}(w_l,\B^{n+1})\subset HF_s^{pq_1}(w_l,\B^{n+1})}$.

\begin{theorem}\label{restriction}
Let $1<p<+\infty$, $1\leq q<+\infty$, $w$ an $A_p$-weight in $\B^{n}$, and $s\in\R$. Then the restriction operator maps $HF_s^{p,q}(w_l,\B^{n+1})$ to $B_{s-\frac1{p}}^{p }(w,\B^{n})$.
\end{theorem}
{\bf Proof of theorem \ref{restriction}:}\par

Since if $q_0\leq q_1\leq +\infty$, ${\displaystyle HF_s^{pq_0}(w_l,\B^{n+1})\subset HF_s^{pq_1}(w_l,\B^{n+1})}$,  it is enough to show that $${\displaystyle HF_s^{p\infty}(w_l,\B^{n+1})_{|\B^{n}}\subset B_{s-\frac1{p}}^{p}(w,\B^{n})}.$$  Let $f\in HF_s^{p\infty}(w_l,\B^{n+1})$, and $k>s$. We have that if $N>0$ is choosen big enough, the representation formula gives that
$$(I+R)^kf(y)=C\int_{\B^{n+1}} (I+R)^kf(z)\frac{(1-|z|^2)^N}{(1-y\overline{z})^{n+2+N}}dv(z).$$
Hence,
$$||f||_{B_{s-\frac1{p}}^{p}(w,\B^{n})}^p 
\preceq \int_{\B^{n}} \left| \int_{\B^{n+1}} (I+R)^kf(z) \frac{(1-|z|^2)^N}{|1-z\overline{y}|^{n+2+N}}dv(z)\right|^p (1-|y|^2)^{(k-s)p}w(y)dv(y).
$$
Next, duality gives that
\begin{equation*}\begin{split}
&||f||_{B_{s-\frac1{p}}^{p}(w,\B^{n})}^p\leq \\&\sup_{ ||\psi||_{L^{p'}(wdv)}\leq 1}\int_{\B^{n}}\int_{\B^{n+1}} |(I+R)^kf(z)|\frac{(1-|z|^2)^N}{|1-z\overline{y}|^{n+2+N}} dv(z)(1-|y|^2)^{k-s} \psi(y)w(y)dv(y)\\&\preceq
\sup_{||\psi||_{L^{p'}(w)}\leq 1}\int_{\B^{n}}\int_{{\bf S}^{n+1}} \int_0^1 |(I+R)^kf(r\zeta)|\frac{(1-r^2)^N}{ |1-r\zeta\overline{y}|^{n+2+N}} dr d\sigma(\zeta)\times\\&(1-|y|^2)^{k-s} \psi(y) w(y)dv(y) \\&\preceq 
 \sup_{||\psi||_{L^{p'}(w)}\leq 1}\int_{{\bf S}^{n+1}}\int_{\B^{n}}
\sup_{r<1} \left( |(I+R)^kf(r\zeta)|(1-r^2)^{k-s}\right)\\& \times \int_0^1 \frac{(1-r^2)^{N+s-k}}{|1-r\zeta\overline{y}|^{ n+2+N}}dr \psi(y) (1-|y|^2)^{k-s}w(y) dv(y) d\sigma(\zeta).
\end{split}\end{equation*}
  
 If $N>0$ is big enough, we have that
${\displaystyle{\int_0^1 \frac{(1-r^2)^{N+s-k}}{|1-r\zeta\overline{y}|^{ n+2+N}}dr\preceq \frac1{|1-\zeta\overline{y}|^{n+1+k-s}}}}$, and the above is bounded by
\begin{equation}\label{firstestimate}
\sup_{ ||\psi||_{L^{p'}(w)}\leq 1}\int_{{\bf S}^{n+1}} \sup_{r<1} \left(    (I+R)^k f(r\zeta)|(1-r^2)^{k-s}\right) \int_{\B^{n}} \frac{\psi(y) (1-|y|^2)^{k-s}}{|1-\zeta\overline{y}|^{n+1+k-s}} w(y) dv(y)d\sigma(\zeta).
\end{equation}

Next, if $M>0$, let $K_M$  be the operator given by  
$$K_M(\psi)(\zeta)=\int_{\B^{n}} \frac{\psi(y) (1-|y|^2)^M}{ |1-\zeta\overline{y}|^{n+1+M}}w(y) dv(y),$$
for $\psi\in L^{p'}(wdv)$, $\zeta\in\overline{\B}^{n+1}$.
We then have that the following lemma holds.
\begin{lemma}\label{estimatebekolle1}
Let $1<p<+\infty$, $M>0$. We then have that $K_M$ is bounded as an operator from $L^{p'}(wdv)$ in $\B^n$ to $L^{p'}({w_l}^{-(p'-1)})$ in ${\bf S}^{n+1}$.
\end{lemma}

Postponing the proof of the lemma , let us finish the proof of the theorem. Applying H\"older's inequality with exponent $p$ to (\ref{firstestimate}),   Lemma   \ref{estimatebekolle1} gives that 
\begin{equation*}\begin{split}
&||f||_{B_{s-\frac1{p}}^{p }(w,\B^{n})}^p\preceq\sup_{||\psi||_{L^{p'}(wdv)}\leq 1}\int_{{\bf S}^{n+1}} \sup_{r<1} \left(  | (I+R)^k f(r\zeta)|(1-r^2)^{k-s}\right)K_{K-s}(\psi)(\zeta)d\sigma(\zeta)\\&\preceq\sup_{ ||\psi||_{L^{p'}(wdv)}\leq 1}\left(\int_{{\bf S}^{n+1}} \sup_{r<1} \left(  | (I+R)^k f(r\zeta)|(1-r^2)^{k-s}\right)^pw(\zeta)d\sigma(\zeta)\right)^\frac1{p}\\&\times
\left(  \int_{{\bf S}^{n+1}}K_{K-s}(\psi)(\zeta)^{p'}w^{-(p'-1)}(\zeta) d\sigma(\zeta)\right)^\frac1{p'} \preceq ||f||_{HF_s^{p\infty}(w,\B^{n+1})}.\qed
\end{split}\end{equation*}

We now give the proof of Lemma \ref{estimatebekolle1}.

{\bf Proof of lemma \ref{estimatebekolle1}:}\par
We observe that for every $\zeta'\in\B^n$, $K_M(\psi)$ is constant on $\Pi^{-1}(\zeta')$ .
Consequently,
\begin{equation*}\begin{split}
&||K_M(\psi)||_{L^{p'}({w_l}^{-(p'-1)}d\sigma)}^{p'}=
\int_{{\bf S}^{n+1}} K_M(\psi)^{p'}(\zeta){w_l}^{-(p'-1)}(\zeta)d\sigma(\zeta)\\&=
\int_{\B^{n}}K_M(\psi)^{p'}(\zeta')w^{-(p'-1)}(\zeta') dv(\zeta')=
||K_M(\psi)||_{L^{p'}(w^{-(p'-1)}dv)}^{p'}.
\end{split}\end{equation*}

Hence, we just have to show that there exists $C>0$ such that for any $\psi\in L^{p'}(wdv)$,
\begin{equation}\label{newestimate}
||K_M(\psi)||_{ L^{p'}(w^{-(p'-1)}dv)}\leq C ||\psi||_{ L^{p'}(wdv)}.
\end{equation}
If we denote by $\psi_1=\psi w$, we have that   $\psi\in  L^{p'}(wdv)$ if and only if $\psi_1\in  L^{p'}(w^{-(p'-1)}dv)$. Thus if $T^*_M$ is the operator defined by
$$\displaystyle T^*_M(\psi_1)(z)= \int_{\B^{n}}\frac{\psi_1(y) (1-|y|^2)^M}{|1-z\overline{y}|^{n+1+M}}dv(y),$$
for $\psi_1 \in L^{p'}(w^{-(p'-1)}dv)$, $z\in \B^{n}$, (\ref{newestimate}) can be rewritten as
\begin{equation}\label{renewestimate}
||T^*_M(\psi_1)||_{ L^{p'}(w^{-(p'-1)}dv)}\leq C ||\psi_1||_{ L^{p'}(w^{-(p'-1)}dv)}.
\end{equation}
And this estimate is again a consequence of proposition 3 in \cite{bekolle}.\qed

Now we prove an extension theorem for weighted holomorphic Besov spaces.

\begin{theorem}\label{extension}
Let $1<p<+\infty$, $s\in\R$, $w$ an $A_p$-weight on $\B^n$. We then have that the extension operator $f\rightarrow f_l$, where  $f_l(z',z_{n+1})=f(z')$, if $(z',z_{n+1})\in{\bf S}^{n+1}$, $z'\in\B^n$, maps boundedly $B_{s-\frac1{p}}^{p}(w,\B^n)$ to $HF_s^{p1}(w_l,\B^{n+1})$.
\end{theorem}
{\bf Proof of theorem \ref{extension}:}\par

Let $f\in B_{s-\frac1{p}}^{p}(w,\B^n)$. If we denote by $R_{n+1}$ the radial derivative in $\B^{n+1}$, then we have
$$(I+R_{n+1})f_l(z',z_{n+1})=(I+\sum_{i=1}^n z_i\frac{\partial}{\partial z_i}+ z_{n+1}\frac{\partial}{\partial z_{n+1}})f_l(z',z_{n+1})= (I+R)f(z'),$$
i.e., $(I+R_{n+1})f_l=\left( (I+R)f\right)_l$.
Consequently, if $k>s$, and $N$ is choosen big enough, duality and the representation formula give
\begin{equation*}\begin{split}&
||f_l||_{HF_s^{p1}(w_l,\B^{n+1})}=\int_{{\bf S}^{n+1}} \left| \int_0^1 \left( (I+R_{n+1})^kf\right)_l(r\zeta)(1-r)^{k-s-1}dr\right|^pw_l(\zeta)d\sigma(\zeta)\\
& \simeq\int_{\B^n}\left|\int_0^1 (I+R)^kf(r\zeta')(1-r)^{k-s-1} dr\right|^p w(\zeta') dv(\zeta')  \\
&=\sup_{ ||\psi||_{L^{p'}(w)}\leq 1}| \int_{\B^n}\left| \int_0^1 (I+R)^k f(r\zeta')(1-r)^{k-s-1}dr\right| \psi(\zeta') w(\zeta') dv(\zeta')|\\&\preceq
\sup_{||\psi||_{L^{p'}(w)}\leq 1} \int_{\B^n}\int_{\B^n}\int_0^1 \frac{|(I+R)^k f(z)|(1-|z|^2)^N(1-r)^{k-s-1} dr}{|1-\overline{z}r\zeta'|^{n+1+N}}dv(z)  \psi(\zeta')w(\zeta')dv(\zeta')\\
& \preceq\sup_{ ||\psi||_{L^{p'}(w)}\preceq 1}  \int_{\B^n}\int_{\B^n}\frac{|(I+R)^k f(z)|(1-|z|^2)^N}{|1-\overline{z}\zeta'|^{n+1+N-(k-s)}}dv(z)\psi(\zeta')w(\zeta')dv(\zeta').
\end{split}\end{equation*}
In order to finish the lemma, an analogous argument to the one used in the restriction theorem, gives that it is enough to show that the mapping defined by
$${\displaystyle{\psi\rightarrow \int_{\B^n} \frac{\psi(\zeta')(1-|z|^2)^{N-(k-s)}}{|1-\overline{z}\zeta'|^{n+1+N-(k-s)}}w(\zeta')dv(\zeta'),}}$$ maps boundedly $L^{p'}(wdV)$ to $L^{p'}(w^{-(p'-1)})$, which is again a consequence of proposition 3 in \cite{bekolle}.\qed

As an immediate consequence of the above two theorems, and the fact that   we have the following relations among the Triebel-Lizorkin spaces \begin{equation*}\begin{split}& HF_s^{p1}(w_l,\B^{n+1})\subset HF_s^{p2}(w_l,\B^{n+1})\subset HF_s^{p\infty}(w_l,\B^{n+1})\\&H_s^{p,2}(w_l,\B^{n+1})=H_s^2(w_l),\end{split}\end{equation*} we obtain the following corollary.
 \begin{corollary}\label{corolariobesov}
 Let $1<p<+\infty$, $s\in \R$, and $w$ a weight in $A_p(\B^n)$. Then the restriction operator from $H_s^p(w_l)$ to $B_{s-\frac1{p}}^{p}(w,\B^n)$ is onto. 
 \end{corollary}
\section{Carleson measures for weighted holomorphic Besov spaces}
In this section we will give a characterization of the Carleson measures for a class of weighted holomorphic Besov spaces.   The main result in \cite{cascanteortega3} shows that for   the Carleson measures for a weighted Hardy-Sobolev space $H_{s_1}^p(w_1)$ for some range of $s_1$ and for a class of weights $w_1$ in ${\bf S}^{n+1}$ coincide with the Carleson measures for the image of the space $L^p(w_1)$  under the operator $K_{s_1}$ of positive kernel given by
$K_{s_1}[L^p[w_1]]$, where
$$K_{s_1}[f](z)=\int_{{\bf S}^{n+1}} \frac{f(\zeta)}{|1-z\overline{\zeta}|^{n+1-s_1}}d\sigma(\zeta).$$  This last problem have been thoroughly studied, (see for instance \cite{sawyerwheedenzhao}).

In our next result, we will see that  we have an analogous situation for Carleson measures for weighted holomorphic Besov spaces.

\begin{theorem}\label{carlesonbesovsobolev}
Let $1<p<+\infty$, $0<s$, $w$ a weight in $B_p(w,\B^n)\cap d_{\tau +1}$,   $0\leq\tau-sp<1$. Let $\mu$ be a  positive Borel measure on $\B^n$. We then have that the following assertions are equivalent:
\begin{itemize}
\item[(i)] There exists $C>0$ such that for any $f\in B_{s}^{p}(w,\B^n)$,
$$||f||_{L^p(d\mu)}\leq C||f||_{B_{s}^{p}(w,\B^n)}.$$

\item[(ii)] There exists $C>0$ such that for any $f\in L^p(wdv)$, 
$$||\int_{\B^n} \frac{f(y)dv(y)}{(1-z\overline{y})^{n+1-(s+\frac1{p})}}||_{L^p(d\mu)}\leq C||f||_{L^p(wdv)}.$$

\item[(iii)] There exists $C>0$ such that for any $f\in L^p(wdv)$, 
$$||\int_{\B^n} \frac{f(y)dv(y)}{|1-z\overline{y}|^{n+1-(s+\frac1{p})}}||_{L^p(d\mu)}\leq C||f||_{L^p(wdv)}.$$

\end{itemize}
\end{theorem}
{\bf Proof of theorem \ref{carlesonbesovsobolev}:}\par
We begin by recalling that  in previous sections we have proved that the regularisation of $w$ is in $A_p(\B^n)\cap D_{\tau+1}$ and that the corresponding Besov spaces coincide. Thus, without loss of generality we may assume that $w$ is in $A_p\cap D_{\tau+1}$. We alsso recall that we have already observed  that if $w\in D_{\tau+1}$, then the lifted weight $w_l$ is in $D_{\tau+1}({\bf S}^{n+1})$.  
We  check that condition (i)  is equivalent to the existence of a constant $C>0$ such that for any $f\in HF_{s+\frac1{p}}^{p,2}(w_l,\B^{n+1})$,
\begin{equation}\label{liftedestimate}
||f||_{L^p(d\mu_l)}\leq C||f||_{HF_{s+\frac1{p}}^{p,2}(w_l,\B^{n+1})},
\end{equation}
where $\mu_l$ is the measure on $\B^{n+1}$ defined by
$\int_{\B^{n+1}}fd\mu_l= \int_{\B^{n}}f(z',0) d\mu(z')$.

Let us show that (i) implies (\ref{liftedestimate}). If $f\in HF_{s+\frac1{p}}^{p}(w_l,\B^{n+1})$, Theorem \ref{restriction} gives that $f_{|B^n}\in B_s^{p}(w,\B^n)$, and that
$||f_{|B^n}||_{B_s^{p}(w,\B^n)}\leq C||f||_{HF_{s+\frac1{p}}^{p,2}(w_l,\B^{n+1})}$.
Since we are assuming that (i) holds, we then have that
$||f_{|B^n}||_{L^p(d\mu)}\leq C||f_{|B^n}||_{B_s^{p}(w,\B^n)}
\leq C||f||_{HF_{s+\frac1{p}}^{p,2}(w_l,\B^{n+1})}$. 
Since $||f_{|B^n}||_{L^p(d\mu)}^p=  
||f||_{L^p(d\mu_l)}^p$, we are done.

Assume now that (\ref{liftedestimate}) holds. Theorem \ref{extension} gives that  $f_l\in HF_{s+\frac1{p}}^{p,2}(w_l,\B^{n+1})$ for any $f\in B_s^{p }(w,\B^n)$, with $||f_l||_{HF_{s+\frac1{p}}^{p,2}(w_l,\B^{n+1})} \leq C|| f||_{B_s^{p }(w,\B^n)}$. The hypothesis on $w$ gives then
$$||f_l||_{L^p(d\mu_l)}\leq C||f_l||_{HF_{s+\frac1{p}}^{p,2}(w_l,\B^{n+1})}\leq C||f||_{B_s^{p }(w,\B^n)}.$$
Since ${\rm sup}\,\, \mu_l\subset \B^n$,
$||f_l||_{L^p(d\mu_l)}^p=||f||_{L^p(d\mu)},$ which gives (i).

Going back to the proof of the theorem, the above observation gives that (i) holds if and only if (\ref{liftedestimate}) does. 
Next, Theorem 2.13 in \cite{cascanteortega3} gives that (\ref{liftedestimate}) can be rewritten as
\begin{equation}\label{holomorphicliftedestimaterealweighted}
||\int_{{\bf S}^{n+1}} \frac{f(\zeta)d\sigma(\zeta)}{(1-z\overline{\zeta})^{n+1-(s+\frac1{p})}}||_{L^p(d\mu_l)}\leq C||f||_{L^p(w_l)}.
\end{equation}

Let us check that (\ref{holomorphicliftedestimaterealweighted}) is equivalent to
\begin{equation}\label{holomorphicestimateweighted}
||\int_{\B^{n}} \frac{f(y)dv(y)}{(1-z\overline{y})^{n+1-(s+\frac1{p})}}||_{L^p(d\mu)}\leq C||f||_{ L^p(w)},
\end{equation}
 for any $f\in  L^p(w)$. Indeed, assume first that (\ref{holomorphicliftedestimaterealweighted}) holds,  let $f\in  L^p(w)$, and let $f_l(z)=f(z')$. We then have that $f_l\in L^p(w_l)$, and $||f_l||_{L^p(w_l)}\simeq ||f||_{ L^p(w)}$.
 Moreover, if $z\in \B^{n}$,
 $$\int_{{\bf S}^{n+1}} \frac{f_l(\zeta) }{(1-z\overline{\zeta})^{n-(s+\frac1{p})}}d\sigma(\zeta)=C
 \int_{\B^n} \frac{f(\zeta')}{(1-z\overline{\zeta'})^{n+1-(s+\frac1{p})}}dv(\zeta'),$$
 and consequently, we obtain (\ref{holomorphicestimateweighted}).

 Assume now that (\ref{holomorphicestimateweighted}) holds, and let $f\in L^p(w_l)$. Then the function $$\widetilde{f}(y)=\frac1{2\pi} \int_{-\pi}^\pi f(y,e^{i\theta}(1-|y|^2))d\theta,$$ for $y\in\B^n$
 is in $ L^p(w)$ (just applying H\"older's inequality) and moreover, $||\widetilde{f}||_{ L^p(w)} \leq C ||f||_{L^p(w_l)}$.
 
In addition,
\begin{equation*}\begin{split}&\int_{\B^n} \frac{\widetilde{f}(y)dv(y)}{(1-z\overline{y})^{n+1-(s+\frac1{p})}}= C\int_{\B^n}\frac{ \int_{-\pi}^\pi f(y,e^{i\theta}(1-|y|^2))d\theta dv(y)}{(1-z\overline{y})^{n+1-(s+\frac1{p})}}=\\&C\int_{{\bf S}^{n+1}}\frac{f(\zeta)}{(1-z\overline{\zeta})^{n+1-(s+\frac1{p})}}d\sigma(\zeta).\end{split}\end{equation*} And that proves that (i) is equivalent to (ii).

Next, the hypothesis   give that $\tau+1-(s+\frac1{p})p=\tau-sp<1$. Since we have observed that the lifted weight $w_l$ satisfies the doubling condition $D_{\tau+1}$, Theorem \ref{weightedtraceinequality}  gives  that  (\ref{liftedestimate}) holds if and only if for any $f\in L^p(w_l)$, $f\geq0$
\begin{equation}\label{liftedestimaterealweighted}
||\int_{{\bf S}^{n+1}} \frac{f(\zeta)d\sigma(\zeta)}{|1-z\overline{\zeta}|^{n+1-(s+\frac1{p})}}||_{L^p(d\mu_l)}\leq C||f||_{L^p(w_l)}.
\end{equation} 
 The same argument used for the holomorphic potential gives that (\ref{liftedestimaterealweighted}) can be rewritten as
 (iii).\qed
 
\begin{corollary}
Let $w_\alpha(z)=(1-|z|)^\alpha$,  $-1<\alpha<p-1$, $s>0$ and $1<p<+\infty$. If    $0< \alpha<p-1$,  assume that $0\leq n+\alpha -sp<1$ and let $\tau = n+\alpha+1$. If $-1<\alpha\leq 0$, assume that $0\leq n-sp<1$ and  let $\tau=n+1$. Let $\mu$ be a  positive Borel measure on $\B^n$. We then have that the following assertions are equivalent:

\begin{itemize}
\item[(i)] There exists $C>0$ such that for any $f\in B_{s}^{p}(w_\alpha,\B^n)$,
$$||f||_{L^p(d\mu)}\leq C||f||_{B_{s}^{p}(w_\alpha,\B^n)}.$$

\item[(ii)] There exists $C>0$ such that for any $f\in L^p(w_\alpha  )$, 
$$||\int_{\B^n} \frac{f(y)dv(y)}{(1-z\overline{y})^{n+1-(s+\frac1{p})}}||_{L^p(d\mu)}\leq C||f||_{L^p(w_\alpha )}.$$

\item[(iii)] There exists $C>0$ such that for any $f\in L^p(w_\alpha  )$, 
$$||\int_{\B^n} \frac{f(y)dv(y)}{|1-z\overline{y}|^{n+1-(s+\frac1{p})}}||_{L^p(d\mu)}\leq C||f||_{L^p(w_\alpha )}.$$

\begin{remark}\label{n=1} We observe that for a wide class of weights in dimension $1$, theorem \ref{carlesonbesovsobolev} holds without imposing any additional condition to the doubling constant $\tau$.  

 {\bf 1.}  Any weight in $A_p({\bf S}^1)$ is in $D_\tau({\bf S}^1)$ for $\tau=p$. Consequentely, lemma \ref{exampleweights} gives that the weight ${\displaystyle\tilde{w}(z)=\frac{1}{(1-|z|) }\int_{I_z}w(\zeta)d\sigma(\zeta)}$ is automatically in $D_{p+1}(\B^1)$. Since $0<\tau-sp=p(1-s)<1$ for $0<s<1$, the hypothesis of theorem \ref{carlesonbesovsobolev} are fullfilled, and a Carleson measure $\mu$ for $B_{s}^p(\tilde{w},\B^1)$ coincides with the ones that satisfy
\begin{equation}\label{casen=1}||\int_{\B^1} \frac{f(y)dv(y)}{|1-z\overline{y}|^{2-(s+\frac1{p})}}||_{L^p(d\mu)}\leq C||f||_{L^p(\tilde{w} dv)}.\end{equation}

{\bf 2.} Let $0<s<1$ and $\max (0,sp-1)\leq \alpha \leq \min (p-1, sp)$. Assume that $\varphi:(0,1]\rightarrow\R$ is a nondecreasing function satisfying that $\varphi(2^kx)\leq C 2^{\alpha k}\varphi(x)$. Let $w_\varphi(z)=\varphi(1-|z|)$. Then proposition \ref{exampledoublingweight}  gives that since $\alpha<p-1$, $w_\varphi\in B_p(\B^1)$, and that for any $\tau>\alpha+2$,  $w_\varphi \in d_\tau(\B^1)$. Thus if we choose $\tau$ such that $ \alpha+2<\tau<2+sp$, $w_\varphi\in B_p(\B^1)\cap d_\tau(\B^1)$. For that choice, $0<\tau-1 -sp<1$, and the agument in remark 4.3.1 can be used to deduce that if $\alpha\in(0,p-1)$, $\mu$ is a Carleson
measure for $B_{s}^p( w_\varphi,\B^1)$ if and only if (\ref{casen=1}) holds replacing $\widetilde{w}$ by $w_\varphi$.
\end{remark}

\end{itemize}\end{corollary}

\end{document}